\documentclass[preprint,12pt]{elsarticle}

\usepackage{amssymb}
\usepackage{amsmath}
\usepackage{geometry}
\usepackage{graphicx}
\usepackage{epstopdf}
\usepackage{algorithm}
\usepackage{algpseudocode}
\usepackage{comment}
\usepackage[utf8]{inputenc}
\usepackage{color}
\usepackage[colorlinks,linkcolor=green,anchorcolor=blue,citecolor=green]{hyperref}

\newtheorem{remark}{Remark}

\def\RR{{\mathbb R}}

\def\cal{\mathcal}

\def\cH{{\cal H}}

\def\cO{{\cal O}}

\def\wt{\widetilde}

\renewcommand{\H}{\mathbf{H}}
\renewcommand{\u}{\mathbf{u}}
\newcommand{\f}{\mathbf{f}}

\newcommand{\x}{\mathbf{x}}

\begin{document}

\begin{frontmatter}



%\title{\textbf{Ray-based Finite Element Method for the Helmholtz Equation in High Frequency Regime}}
%\title{\textbf{Learning adaptive basis for ray-based Finite Element Method for the High-Frequency Helmholtz Equation}}
\title{\textbf{Learning Dominant Wave Directions For Plane Wave Methods For High-Frequency Helmholtz Equations}}

\author{Jun Fang}
\address{Department of Mathematics, University of California, Irvine}

\author{Jianliang Qian}
\address{Department of Mathematics, Michigan State University}

\author{Leonardo Zepeda-N\'{u}\~{n}ez}
\address{Department of Mathematics, University of California, Irvine}

\author{Hongkai Zhao}
\address{Department of Mathematics, University of California, Irvine}

\begin{abstract}
% We present a ray-based finite element method (Ray-FEM) for the high-frequency Helmholtz equation with smooth wave-speed.
% The key idea is to probe the medium with the same source but using a relative low frequency wave to extract/learn the ray directions by applying Numerical Micro-Local Analysis (NMLA) to the computed wave field. Then the ray information is incorporated into the Galerkin space to approximate the high frequency Helmholtz solution. Moreover, the process can be continued to further improve the approximation for both the ray directions and the high frequency wave field iteratively. Our method only requires the minimum degrees of freedom, i.e., a fixed number of grid points per wave length, to achieve both stability and expected accuracy for the high frequency Helmholtz solution without the usual pollution effect. Fast solvers are developed to  solve the resulting linear systems to achieve the overall efficiency of our method. Numerical tests in both 2D and 3D are presented to corroborate our method.

We present a ray-based finite element method (ray-FEM) by learning basis adaptive to the underlying high-frequency Helmholtz equation in smooth media. Based on the geometric optics ansatz of the wave field, we learn local dominant ray directions by probing the medium using low-frequency waves with the same source. Once local ray directions are extracted, they are incorporated into the finite element basis to solve the high-frequency Helmholtz equation. This process can be continued to further improve approximations for both local ray directions and the high frequency wave field iteratively.
The method requires a fixed number of grid points per wavelength to represent the wave field and achieves an asymptotic convergence
as the frequency $\omega\rightarrow \infty$ without the pollution effect.
A fast solver is developed for the resulting linear system with an empirical complexity  $\mathcal{O}(\omega^d)$ up to a poly-logarithmic factor. Numerical examples in 2D are presented to corroborate the claims.
\end{abstract}


\end{frontmatter}



 %!TEX root = draft.tex

\section{Introduction}
We consider the Helmholtz equation:
\begin{equation}\label{eq:Helmholtz}
\cH u  := -\Delta u(\textbf{x}) - \frac{\omega^2}{c^2(\textbf{x})}u(\textbf{x}) = f(\textbf{x}), \quad \x \in \Omega\subseteq \RR^d,
\end{equation}
plus boundary (or radiation) conditions, where $\omega$ is the frequency, $c(\textbf{x})>0$ is the wave speed, and $f(\textbf{x})$ is the source distribution.

The numerical solution of the Helmholtz equation \eqref{eq:Helmholtz} in the high-frequency regime, i.e., $\omega \gg 1$, is notoriously hard to compute. From Shannon's sampling principle \cite{Shanon:Communication_In_The_Presence_Of_Noise}, to resolve a general wave field oscillating at frequency $\omega$, a mesh size $h=\cO(\omega^{-1})$ is necessary and sufficient. Hence the intrinsic degrees of freedom is $\cO(\omega^{d})$, implying that the theoretical optimal overall complexity to solve \eqref{eq:Helmholtz} is $\cO(\omega^{d})$. In general, an overall complexity of optimal order is difficult to achieve in practice due to two typical challenges: how to design a discretization method that can achieve both accuracy and stability with the minimum order of degrees of freedom; and how to solve the linear system issued from the discretization with an almost linear complexity as the frequency becomes large.

% The numerical solution of the Helmholtz equation \eqref{eq:Helmholtz} in the high-frequency regime, i.e., $\omega \gg 1$, is notoriously hard to compute. There exists two typical challenges: how to design a discretization method that can achieve both accuracy and stability with the minimum degrees of freedom (shall we explain what are the minimum dof?); and how to solve the linear system issued from the discretization efficiently.
% From Shannon's sampling principle \cite{Shanon:Communication_In_The_Presence_Of_Noise}, to resolve a general wave-field which is band-limited at $\omega$, a mesh size $h=\cO(\omega^{-1})$ is necessary and sufficient. Hence the theoretical optimal overall complexity to solve \eqref{eq:Helmholtz} is $\cO(\omega^{d})$. In general, the optimal overall complexity is difficult to achieve in practice due to the above two challenges.

Numerically, it is difficult to design a sparse discretization that can achieve both accuracy and stability under the condition $\omega h = \cO(1)$ as $\omega$ becomes large. This is usually called pollution effect in error estimates for finite element methods \cite{Babuska:Is_Pollution_Effect_Aviodable,Babuska:A_Generalized_Finite_Element_Method_for_solving_the_Helmholtz_equation_in_two_dimensions_with_minimal_pollution,Ihlenburg_Babuska:Solution_of_Helmholtz_problems_by_knowledge_based_FEM}, i.e., the ratio between numerical error and best approximation error from a discrete finite element space is $\omega$ dependent.

A popular strategy to solve the high-frequency Helmholtz equation \eqref{eq:Helmholtz} is to incorporate appropriate oscillatory behaviors into the basis for Galerkin methods in order to represent the oscillatory wave field governed by the Helmholtz equation more accurately and efficiently. The methods based on this strategy are usually called wave based methods. The key issue for this strategy is how to design the oscillatory basis. Several approaches have been developed along this direction. One approach is to incorporate plane waves with a predetermined even distribution in directions into the basis. For example, products of plane waves with local finite element basis are used in the  generalized finite element method \cite{Babuska:A_Generalized_Finite_Element_Method_for_solving_the_Helmholtz_equation_in_two_dimensions_with_minimal_pollution, Strouboulis_Babuska_Hidajat_06:GFEM}. Trefftz-type methods use local solutions of the Helmholtz equation \cite{Hiptmair_Moiola_Perugia:A_Survey_of_Trefftz_Methods_for_the_Helmholtz_Equation} as the basis functions, which in the case of piece-wise constant media are plane waves.

A challenging issue in these approaches is how to choose the number of plane wave directions {\it\`{a} priori}. It is well known that in order to achieve a good accuracy, a fine, $\omega$ dependent, resolution in angle space is required. However, since many of these directions may not be relevant to a specific problem setup, it will not only increase the degrees of freedom significantly but also make the resulting linear system ill-conditioned due to the numerical loss of orthogonality of the elements of the basis. %, due to numerical redundancy.
For heterogeneous media, the construction of local solutions of the Helmholtz equation becomes more challenging, and some recent papers have focused on how to construct such local solution: a generalized plane wave method was developed in \cite{Imbert_Monk:Numerical_simulation_of_wave_propagation_in_inhomogeneous_media_using_Generalized_Plane_Waves}; and a two-scale method using local numerical solutions of the Helmholtz equation on a fine mesh as the basis function on a coarse mesh was developed in \cite{Gallistl_Peterseim:Stable_Multiscale_Petrov-Galerkin_Finite_Element_Method_for_High_Frequency_Acoustic_Scattering}.

Another approach to solve the high-frequency Helmholtz equation is based on the geometric optics ansatz of the wave field. In the ansatz, the high-frequency wave field is approximated by a superposition of a finite number, that we denote by $N$, of dominant wave fronts of the form
\begin{equation}
    u(\textbf{x}) \approx \sum_{n=1}^N A_n(\textbf{x}) e^{i\omega\phi_n(\textbf{x})}
\end{equation}
at each point; see \cite{Luo:Fast_Huygens_sweeping_methods_for_Helmholtz_equations_in_inhomogeneous_media_in_the_high_frequency_regime} for examples. The main advantage is that both the amplitude functions $A_n(\x)$ and phase functions $\phi_n(\x)$ are independent of the frequency $\omega$ and hence are non-oscillatory and smooth except at a measure zero set, e.g., focus points, caustics, corners in a smooth medium. However, both the amplitude and phase functions depend on the medium, source distribution and domain boundary, i.e., they are global and problem specific. Once the phase functions of the wave fronts are available, the oscillatory pattern of the wave field is known and can be incorporated into the numerical approximation to improve both stability and accuracy. Instead of incorporating many problem independent plane waves to the basis, the main idea of phase based numerical methods is to explicitly incorporate the phase function into the basis functions \cite{Giladi07:Asymptotically_derived_boundary_elements_for_the_Helmholtz_equation_in_high_frequencies, DBLP:journals/jcphy/NguyenPRC15}.

The main issue in phase-based methods is how to compute the phase functions, which are global unknown functions depending on the medium and the source distribution. Only in simple geometry and homogenous medium the phase function can be computed analytically; otherwise, it has to be computed numerically \cite{Luo:Fast_Huygens_sweeping_methods_for_Helmholtz_equations_in_inhomogeneous_media_in_the_high_frequency_regime}. Usual approaches include Lagrangian methods, such as ray tracing by solving ordinary differential equations (ODE), and Eulerian methods based on partial differential equation (PDE), e.g., the Eikonal equation, satisfied by the phase function. Both these methods face numerical difficulties for general media in practice. For example, convergence and divergence of rays cause difficulties to obtain an accurate and uniform ray field. Since the solution to the PDE is single valued, it can only capture the phase function corresponding to the first arrival time of the wave front. For a general medium with varying speed, computing the appropriate global phase functions $\phi_n(\x)$ is a challenging task since different phase functions may be defined in different regions the boundaries of which are difficult to determine. In other words, the number of terms in the geometric optics ansatz varies from point to point; see \cite{Luo:Fast_Huygens_sweeping_methods_for_Helmholtz_equations_in_inhomogeneous_media_in_the_high_frequency_regime} for examples. Moreover, since the phase functions are precomputed once for all and fed into the computation of the solution for the Helmholtz equation in these approaches, the phase function has to be computed extremely accurate because the numerical error is amplified by $\omega$. Instead of incorporating the global phase function into the basis, one can use local ray directions, $\widehat{\textbf{d}}_n(\x) = \frac{\nabla \phi_n(\x)}{|\nabla \phi_n(\x)|}$, as local dominant plane wave direction of $e^{i\omega \widehat{\textbf{d}}_n \cdot \x}$ and incorporate them into the local basis. Compared to the usual plane wave methods that typically require a large number of evenly distributed wave directions independent of the problem, only directions of dominant wave fronts relevant to the problem are involved. Hence the degrees of freedom can be kept minimal and ill-conditioning of the resulting linear system due to redundancy can be reduced.
This combines the advantages of plane wave methods and phase-based methods. As is shown in \cite{Betcke_Phillips:Approximation_by_dominant_wave_directions_in_plane_wave_methods}, using dominant wave directions in plane wave methods can significantly improve efficiency and accuracy for solving the high-frequency Helmholtz equation in heterogeneous medium. Moreover, it was shown that there is some tolerance in the approximation of the dominant wave directions. For those examples in \cite{Betcke_Phillips:Approximation_by_dominant_wave_directions_in_plane_wave_methods}, a simple ray tracing method was used to find the rays and the corresponding dominant wave directions at each point for a very simple setup in a homogeneous medium and the exact phase function was given for a heterogeneous medium.

Although incorporating problem specific dominant wave directions into basis functions has shown very promising results, a key issue in practice is  how to find the dominant wave directions in a stable, efficient and systematic way for general media.

In this work we propose an approach that is based on learning the dominant wave directions specific to the medium and source distribution. In particular, we probe the same medium by the same source, i.e., solving the Helmholtz equation \eqref{eq:Helmholtz} with the same $c(\x), f(\x)$, for a relative low frequency $\widetilde{\omega} \sim \sqrt{\omega}$. The computed wave field is post-processed by numerical micro-local analysis (NMLA) or other signal processing tools to estimate dominant wave directions. These estimated wave directions are then used in plane wave methods to solve the original high-frequency Helmholtz equation. In our approach, global phase functions are not needed. Instead, dominant wave directions are computed locally where both the number of dominant wave directions and the directions can vary from point to point. It provides the flexibility to deal with general media. Moreover, once a more accurate wave field is computed, it can be used to get a better estimation of the dominant wave directions and then used to improve the high-frequency wave field again. In other words, both the dominant wave directions and the high-frequency wave field can be computed and improved in an iterative way. In this paper, we develop a simple ray-based finite element (ray-FEM) method in 2D for smooth media as a proof of concept study of our proposed approach. We start with a finite element mesh with mesh size $h$ satisfying $wh=\cO(1)$, i.e., a few points per wavelength. First, the low frequency Helmholtz equation with $\widetilde{\omega} \sim \sqrt{\omega}$ is solved on the mesh with quasi-optimality since $\widetilde{\omega}^2h=\cO(1)$. Then NMLA \cite{Benamou:NMLA_harmonic_wavefields, Benamou:NMLA_revisited} (see Section \ref{section:NMLA}) is applied to the computed low-frequency wave field to estimate local dominant wave directions. In the second step we incorporate these estimated dominant wave directions into the local finite element basis as the generalized finite element method to solve the high-frequency Helmholtz equation on the same mesh. If necessary, the solution for the high-frequency Helmholtz equation can also be processed by NMLA to improve the estimate of local dominant wave directions which can be used to further improve the high-frequency solution. We also develop an iterative fast solver for the discretized linear system based on polarized traces which can achieve $\cO(N)$ complexity with a possible poly-logarithmic factor for smooth media, where $N$ is the total number of unknowns. We use numerical examples to show that our approach can achieve the following goals:
(1) a stable discretization with degrees of freedom of the minimal order $\cO(\omega^d)$, (2) a fast solver with approximate linear complexity for the discretized linear systems, (3) an asymptotic convergence rate of $\cO(\omega^{ -\frac{1}{2}})$ as $\omega\rightarrow \infty$.

Here is an outline of this paper. We first describe the ray-FEM using geometric optics ansatz as the motivation and study its approximation property in Section \ref{section:ray_FEM}. In Section \ref{section:NMLA} we introduce the NMLA with its stability and local ray direction error analyzed in Appendix A and B. Section \ref{Section:Algorithms} provides the full presentation of the numerical algorithm whose complexity is given in Section \ref{section:complexity}. Numerical results are presented in Section \ref{Section:Numerical_Experiments}. Conclusion and future works are summarized in Section \ref{Section:Conclusion}.

\section{The Ray-FEM Method} \label{section:ray_FEM}

In this section we describe the ray-FEM method for the Helmholtz equation and its rationale based on the geometric optics ansatz. We explain briefly the ansatz and approximate it locally via a superposition of plane waves propagating in dominant directions. These plane waves with their associated dominant directions are incorporated into the finite element basis functions to improve both stability and accuracy for the computation of the solution to the high-frequency Helmholtz equation.

In this section we suppose that the dominant directions are known exactly. However, in Section \ref{section:NMLA} we will describe how to learn the dominant wave directions by probing the medium using low-frequency waves.

 % In next section we will describe how to learn the dominant wave directions by probing the medium using low frequency waves. Once the estimation of dominant wave directions is available, plane waves with these dominant directions will be incorporated into the finite element basis functions to improve both stability and accuracy for the computation of high frequency Helmholtz equation.

We use the following boundary value problem in 2D for the illustration of our method,
\begin{equation} \label{eq:Helmtholtz_model}
\left\{ \begin{array}{rl}
  -\Delta u - k^2(\x) u = f,  						 &  \mbox{in} \quad \Omega,  \\
\frac{\partial u}{\partial n} + i\beta k(\x) u  = g, &  \mbox{on} \quad \partial \Omega,
  \end{array} \right.
\end{equation}
where $\Omega$ is a bounded Lipschitz domain in $\mathbb{R}^2$, and $k(\x) = \omega/c(\x)$ is the inhomogeneous wave number. \eqref{eq:Helmtholtz_model} is usually  refered to as the Helmholtz equation with impedance boundary conditions. This equation was chosen in order to easily impose another types of boundary conditions by modifying the coefficient $\beta$. Specifically, the Dirichlet boundary condition corresponds to $\beta = \infty$ and the first order absorbing boundary condition to $\beta = \pm 1$. Moreover, it is easy to extend \eqref{eq:Helmtholtz_model} to incorporate absorbing boundary conditions implemented via PML \cite{Berenger:PML}, as it will be performed in the numerical experiments in Section \ref{section:experiments_fast_methods}.

%We first introduce the geometrical optics Ansatz , which constitutes the backbone of the algorithm. Then we approximate the geometrical optics Ansatz locally via a superposition of plane waves, which render the approximation suitable to enrich a low order finite element space, which is properly defined afterwards. Finally, we show some basic approximation results using the enriched finite element space.
%
%We point out that in this section we suppose that the phase information is {\it known}, and we use it to build the approximation spaces. In the sequel, mostly in Section \ref{section:NMLA}, we will describe how to obtain the phase information.

\subsection{Geometric optics ansatz and local plane wave approximation } \label{section:geometrical_optics}
The standard derivation of geometric optics ansatz is the use of WKJB approximation \cite{Rayleigh:On_the_Propagation_of_Waves_through_a_Stratified_Medium_with_Special_Reference_to_the_Question_of_Reflection,Jeffreys:On_Certain_Approximate_Solutions_of_Linear_Differential_Equations_of_the_Second_Order}  (or the L\"uneberg-Kline expansion \cite{Kline_Kay:Electromagnetic_Theory_and_Geometrical_Optics}) for the solution to the Helmholtz equation \eqref{eq:Helmholtz}:
\begin{equation} \label{eq:WKJB}
	u(\x) \sim  e^{i\omega \phi (\x)} \sum_{\ell = 0}^{\infty} \frac{A_{\ell}(\x)}{\omega^{\ell} }.
\end{equation}
By taking $\omega \rightarrow \infty $ and considering only the first term one has
\begin{equation}
u(\textbf{x}) = A(\textbf{x}) e^{i\omega\phi(\textbf{x})}  + O \left( \frac{1}{\omega} \right),
\end{equation}
where $A$ is usually called the amplitude and $\phi$ the phase. The key features of the geometric optics ansatz are:
\begin{itemize}
\item $A$ and $\phi$ are {\it independent} of the frequency $\omega$;
\item $A$ and $\phi$ {\it depend} on the medium, $c(\x)$; and the source distribution, $f(\x)$.
\end{itemize}
Moreover, except for a small set of points, e.g.,  source/focus points, caustics, and discontinuities of the medium, they are smooth functions satisfying the following PDE system:
\begin{equation}
\mbox{(eikonal)} \quad |\nabla \phi | = \frac{1}{c},   \qquad \mbox{(transport)} \quad 2\nabla\phi \cdot \nabla A + A \Delta \phi = 0.
\label{eq:eikonal}
\end{equation}

In general, the phase function, $\phi$, and the amplitude function, $A$, are multi-valued functions corresponding to multiple arrivals of wave fronts. Hence one can further decompose the geometric optics ansatz into a superposition of a number of wave fronts in the form:
\begin{equation} \label{eq:plane_wave_expansion}
u(\textbf{x}) = \sum_{n=1 }^{N(\x)} A_{n}(\textbf{x}) e^{i\omega\phi_{n}(\textbf{x})}  + O \left( \frac{1}{\omega} \right),
\end{equation}
where $N(\x)$ is the number of fronts/rays passing through $\x$, and the phases $\phi_n$ and amplitudes $A_n$ are single valued functions satisfying the eikonal/transport equations (\ref{eq:eikonal}), each defined in a suitable domain with suitable boundary conditions \cite{Benamou:An_Introduction_to_Eulerian_Geometrical_Optics}.

Based on the above geometric optics ansatz, one can derive a local plane wave approximation at any point where $\phi_n$ and $A_n$ are smooth with variations on a $\cO(1)$ scale. Indeed, using Taylor expansions on a small neighborhood around an observation point $\textbf{x}_0$, we have
%\begin{align}
%A_{n}(\textbf{x}) 		& =    A_{n}(\mathbf{x}_{0}) + \nabla A_{n}(\mathbf{x}_{0}) \cdot (\x - \x_0) + \mathcal{O} (h^2) ,\nonumber \\  \label{eq2.4}
%\phi_{n}(\mathbf{x}) 	& = \phi_{n}(\mathbf{x}_{0}) + \nabla \phi (\mathbf{x}_{0}) \cdot (\mathbf{x}-\mathbf{x}_{0}) +  \mathcal{O}(h^2) \nonumber.
%\end{align}
%
%Which provides the following approximation
\begin{equation} \label{eq:plane_wave_approximation_with_nabla}
u(\x) \!\!=\!\! \!\sum_{n=1 }^{N(\x_0)} \!\!\left ( A_{n}(\textbf{x}_{0})\! +\! \nabla \!A_{n}(\textbf{x}_{0}) (\x \!-\! \x_0) \right )  e^{ i\omega \left (  \phi_{n} (\x_0) + \nabla \phi(\x_0) \cdot (\x - \x_0 ) \right )  }\!  +\! \mathcal{O} \left(\! h^2 \!+\! \omega h^2 \!+\! \frac{1}{\omega}\! \right),
\end{equation}
for $|\x-\x_0|<h\ll 1$.

%We denote by  $\textbf{d}_{n} \in \mathbb{S}^1 = \lbrace (\cos \theta, \sin \theta): \theta  \in [0, 2\pi) \rbrace$ the direction of propagation of the rays at $\textbf{x}_0$ (or dominant ray directions), i.e.,  $\nabla \phi_{n}(\textbf{x}_0) =\frac{\textbf{d}_{n}}{c(\textbf{x}_0)}$; and $k(\x_0) = \omega/c(\x_0)$. Using the notations just introduced we can rewrite the approximation in \eqref{eq:plane_wave_approximation_with_nabla} as a local superposition of plane waves, which results in
Define
\begin{equation} \label{eq:dominant_direction}
   \widehat{\textbf{d}}_{n}:=\frac{\nabla \phi_n(\x_0)}{|\nabla \phi_n(\x_0) |} = c(\x_0) \nabla \phi_n(\x_0)
\end{equation}
as the ray directions of the wave fronts at $\x_0$, $k(\x_0) = \omega/c(\x_0)$, and
\begin{equation} \label{eq:amplitude}
    B_{n}(\x) = (A_{n}(\textbf{x}_{0}) + \nabla A_{n}(\textbf{x}_{0}) (\x - \x_0) )e^{i\omega (\phi_{n}(\textbf{x}_{0})-\nabla \phi(\x_0)\cdot \x_0)}
\end{equation}
the affine complex amplitude. By replacing \eqref{eq:dominant_direction} and \eqref{eq:amplitude} in \eqref{eq:plane_wave_approximation_with_nabla} we have
\begin{equation} \label{eq:plane_wave_approximation}
u(\x) = \sum_{n=1}^{N(\x_0)} B_{n}(\x) e^{ i k(\x_0)  \widehat{\textbf{d}}_{n} \cdot  \x }  + \mathcal{O} \left( h^2 + \omega h^2 + \frac{1}{\omega} \right),
\end{equation}
for $|\x-\x_0|<h\ll 1$. From \eqref{eq:plane_wave_approximation} we have that $u$ can be approximated locally by a superposition of plane waves propagating in certain directions with affine complex amplitudes. Moreover, as $\omega \rightarrow \infty$, such that $ \omega h=\cO(1)$, the asymptotic error for the local plane wave approximation \eqref{eq:plane_wave_approximation} is $\cO(\omega^{-1})$, which is of the same order as the asymptotic error for the original geometric ansatz \eqref{eq:plane_wave_expansion}. We use \eqref{eq:plane_wave_approximation} as the motivation to construct local finite element basis with mesh size $h=\cO(\omega^{-1})$, in which an affine function is multiplied by plane waves oscillating in those ray directions, resulting in local approximations similar to \eqref{eq:plane_wave_approximation}.

\subsection{Ray-based FEM formulation}
We use a finite element method to compute the solution to \eqref{eq:Helmtholtz_model} whose standard weak formulation is given by
\begin{equation} \label{eq:Helmholtz_equation_weak}
\mbox{Find } u \in H^1(\Omega), \mbox{ such that  }  \mathcal{B}(u,v) = \mathcal{F}(v), \quad \forall v \in H^1(\Omega),
\end{equation}
where
\begin{equation}
\mathcal{B}(u,v) := \int_{\Omega}\nabla u \cdot \nabla \overline{v} dV - \int_{\Omega} k^2u\overline{v}dV + i\beta \oint_{\partial \Omega} k u \overline{v}dS,
\end{equation}
\begin{equation}
\mathcal{F}(v) := \int_{\Omega}f\overline{v}dV + \oint_{\partial \Omega}g\overline{v}dS.
\end{equation}

The domain, $\Omega$, is discretized with a standard regular triangulated mesh, with mesh size $h$, which we denote by $\mathcal{T}_h = \{K\}$, where $K$ represents a triangle of the mesh. Using the aforementioned mesh we define two approximation spaces for the variational formulation \eqref{eq:Helmholtz_equation_weak}:

\begin{itemize}
\item Standard FEM (S-FEM), we use low order $\mathbb{P}1$ finite elements, i.e., piece-wise bilinear functions,
\item Ray-FEM, we use $\mathbb{P}1$ finite elements multiplied by plane waves as in  \eqref{eq:plane_wave_approximation}.
\end{itemize}

For a given element $K \in \mathcal{T}_h$, we denote by $V_j$ and $\textbf{x}_j, j = 1,2,3$, the vertices of $K$ and their coordinates respectively. Moreover, we denote by $\{\varphi_j(\x) \}_{j = 1}^3$ a partition of unity consisting of piecewise bilinear functions satisfying  $\varphi_j(\x_i) = \delta_{ij}$, $i, j = 1,2,3$, where $\delta_{ij}$ is the Kronecker delta. The basis given by $\{\varphi_j(\x) \}_{j = 1}^3$ is usually called the nodal basis for Lagrange $\mathbb{P}1$ finite elements.
The standard local approximation space is given by
\begin{equation}
	V_{S}(K) = \mbox{span} \{ \varphi_j(\textbf{x}),  j = 1,2,3  \},
\end{equation}
and the global $\mathbb{P}1$ finite element space
\begin{equation}
	V_{S}(\mathcal{T}_h) =  \{ v\in C^0(\overline{\Omega}) : v|_K \in V_S(K), \forall K \in \mathcal{T}_h \}.
\end{equation}

To define the ray-FEM we enrich the $\mathbb{P}1$ finite elements by incorporating the ray information. Let $\{ \widehat{\textbf{d}}_{j,l}\}_{l = 1}^{n_j}$ be $n_j$ ray directions at the vertex $V_j$, define the ray-based local approximation space by
\begin{displaymath}
V_{Ray}(K) = \mbox{span} \{ \varphi_j(\textbf{x})e^{ik_j\widehat{\textbf{d}}_{j,l}\cdot \textbf{x} }, \quad k_j = k(\textbf{x}_j), \quad j = 1,2,3, \quad l = 1,...,n_j  \},
\end{displaymath}
and the global ray-FEM space by
\begin{displaymath}
V_{Ray}(\mathcal{T}_h) =  \{ v\in C^0(\overline{\Omega}) : v|_K \in V_{Ray}(K), \forall K \in \mathcal{T}_h \}.
\end{displaymath}

We can define the Standard FEM method by
\begin{equation} \label{S-FEM}
\mbox{Find } u \in V_{S}(\mathcal{T}_h), \mbox{ such that  }  \mathcal{B}(u,v) = \mathcal{F}(v), \quad \forall v \in V_{S}(\mathcal{T}_h).
\end{equation}
Analogously, we define the ray-FEM method by
\begin{equation} \label{R-FEM}
\mbox{Find } u \in V_{Ray}(\mathcal{T}_h), \mbox{ such that  }  \mathcal{B}(u,v) = \mathcal{F}(v), \quad \forall v \in V_{Ray}(\mathcal{T}_h).
\end{equation}

\subsection{Approximation property of ray-FEM with exact ray information} \label{section:approximation_L2_norm}

We provide a simple computation to estimate the approximation error of the ray-FEM space. In particular, we compute an asymptotic bound on $\inf_{u_h \in V_{Ray}(\mathcal{T}_h)} || u - u_h ||_{L^2(\Omega)}$, which is achieved by estimating the interpolating error using $V_{Ray}(K)$ as a basis.

In the computation we assume that the ray direction, which is the gradient of the phase function $\phi$, and the phase function itself, are exactly known. For simplicity, we assume $N = 1$ for the asymptotic formula in \eqref{eq:plane_wave_expansion}, i.e., that only one ray crosses each point of the domain. In addition we suppose that $f$, the source, is zero inside the domain. Under those circumstances $A$ and $\phi$ are smooth; given that the source is outside the domain there is no singularity in the amplitude, and given that only one ray crosses each point in the domain, no caustic occurs. From the geometric optics ansatz, we have
\begin{equation}
  u(\textbf{x}) =  A(\textbf{x}) e^{i\omega\phi(\textbf{x})}  + O \left( \omega^{-1} \right).
\end{equation}

Let $\Omega$ be our domain of interest, $\mathcal{T}_h$ a triangle mesh of $\Omega$ with width $h$. We denote by $N_h$ the total number of vertices on the mesh $\mathcal{T}_h$, $ \lbrace \x_j \rbrace_{j = 1}^{N_h}$ and $\lbrace \varphi_j(\x) \rbrace_{j = 1}^{N_h}$ are the coordinates of all mesh nodes and their corresponding nodal basis functions for standard $\mathbb{P}1$ element.

We note that $ e^{i\omega [ \phi(\textbf{x}_j) - \nabla \phi(\textbf{x}_j) \cdot \textbf{x}_j ]}$ is a constant for the nodal basis associated to $\x_j$ in an element $K$. From this observation we can easily deduce that the local ray-FEM space can be rewritten as
\begin{displaymath}
\begin{array}{ll}
V_{Ray}(K) & =  \mbox{span} \{ \varphi_j(\textbf{x})e^{ik_j\textbf{d}_{j}\cdot \textbf{x} } \} =  \mbox{span} \{ \varphi_j(\textbf{x})e^{i\omega \nabla \phi (\textbf{x}_{j}) \cdot \textbf{x} } \}  \\[1.0ex]
&=\mbox{span} \{ \varphi_j(\textbf{x})e^{i\omega \nabla \phi (\textbf{x}_{j}) \cdot \textbf{x} } e^{i\omega [ \phi(\textbf{x}_j) - \nabla \phi(\textbf{x}_j) \cdot \textbf{x}_j ]} \} \\[1.0ex]
& =  \mbox{span} \{ \varphi_j(\textbf{x}) e^{i\omega [ \phi(\textbf{x}_j) +  \nabla \phi(\textbf{x}_j) \cdot ( \textbf{x} - \textbf{x}_j) ]} \}.
\end{array}
\end{displaymath}
Hence the nodal interpolation of the solution can be written as
\begin{equation} \label{nodal_interpolation}
u_I = \sum_{j=1 }^{N_h} A(\textbf{x}_j)\varphi_j(\textbf{x}) e^{i\omega [ \phi(\textbf{x}_j) +  \nabla \phi(\textbf{x}_j) \cdot ( \textbf{x} - \textbf{x}_j) ]},
\end{equation}
which, by construction lies within the global ray-FEM space $V_{Ray}(\mathcal{T}_h)$.

Let $S_j$ be the support of $\varphi_j(\textbf{x})$, and $|S_j| \sim O(h^2)$ be the area of $S_j$. Then using the triangular inequality and the smoothness assumptions we have
\begin{displaymath}
\begin{array}{ll}
\Vert u - u_I \Vert_{L^2(\Omega)} & \leq \Vert A(\textbf{x})e^{i\omega \phi (\textbf{x})} -  \sum_{j=1 }^{N_h} A(\textbf{x}_j)\varphi_j(\textbf{x}) e^{i\omega  \phi(\textbf{x})} \Vert_{L^2(\Omega)} \\[1.0ex]
&\quad  + \Vert  \sum_{j=1 }^{N_h} A(\textbf{x}_j)\varphi_j(\textbf{x}) \left( e^{i\omega  \phi(\textbf{x})}-   e^{i\omega [ \phi(\textbf{x}_j) +  \nabla \phi(\textbf{x}_j) \cdot ( \textbf{x} - \textbf{x}_j) ]}   \right) \Vert_{L^2(\Omega)} + O(\omega^{-1})\\[1.5ex]
& \leq  \Vert A(\textbf{x}) -  \sum_{j=1 }^{N_h} A(\textbf{x}_j)\varphi_j(\textbf{x}) \Vert_{L^2(\Omega)} \\[1.0ex]
& \quad +   \sum_{j=1 }^{N_h} \Vert A  \Vert_{L^{\infty}(\Omega)} \Vert e^{i\omega  \phi(\textbf{x})}-   e^{i\omega [ \phi(\textbf{x}_j) +  \nabla \phi(\textbf{x}_j) \cdot ( \textbf{x} - \textbf{x}_j) ]}    \Vert_{L^2(S_j)} + O(\omega^{-1})
\\[1.5ex]
& \lesssim h^2 \vert A \vert_{H^2(\Omega)} +  \sum_{j=1 }^{N_h} \Vert A  \Vert_{L^{\infty}(\Omega)} \omega h^2 \Vert \nabla^2 \phi  \Vert_{L^{\infty}(\Omega)} |S_j| + O(\omega^{-1}) \\[1.5ex]
& \lesssim h^2 \vert A \vert_{H^2(\Omega)} +  \omega h^2 \Vert A  \Vert_{L^{\infty}(\Omega)}  \Vert \nabla^2 \phi  \Vert_{L^{\infty}(\Omega)}  +  O(\omega^{-1}). \\
\end{array}
\end{displaymath}
This implies
\begin{equation} \label{approximation_property}
\inf_{u_h \in V_{Ray}(\mathcal{T}_h)} || u - u_h ||_{L^2(\Omega)} \lesssim h^2 \vert A \vert_{H^2(\Omega)} +  \omega h^2 \Vert A  \Vert_{L^{\infty}(\Omega)}  \Vert \nabla^2 \phi  \Vert_{L^{\infty}(\Omega)}  +  O(\omega^{-1}) .
\end{equation}
Or, asymptotically,
\begin{equation}  \label{eq:approximation_property_exact_ray}
    \inf_{u_h \in V_{Ray}(\mathcal{T}_h)} || u - u_h ||_{L^2(\Omega)} \lesssim \cO(h^2 + \omega h^2 + \omega^{-1} ).
\end{equation}
If the exact rays are known and the mesh size follows $h \sim \omega^{-1}$, then we have
\begin{equation}
    \inf_{u_h \in V_{Ray}(\mathcal{T}_h)} || u - u_h ||_{L^2(\Omega)} \lesssim \cO(\omega^{-1}),
\end{equation}
i.e. that the approximation error decays linearly with $\frac{1}{\omega}$, without oversampling.

\begin{remark}
The ray information can be incorporated into other Galerkin basis in the same fashion. For example, in the hybrid numerical asymptotic method of \cite{Giladi_Keller:A_hybrid_numerical_asymptotic_method_for_scattering_problems}, the basis functions are constructed by multiplying nodal piece-wise bilinear function to oscillating functions with phase factors; the plane wave DG method of  \cite{Betcke_Phillips:Approximation_by_dominant_wave_directions_in_plane_wave_methods} employs the products of small degree polynomials and dominant plane waves as basis functions; the phase-based hybridizable DG method of \cite{DBLP:journals/jcphy/NguyenPRC15} considers basis functions as products of polynomials and phase-based oscillating functions. However, the phase or ray information in these methods is obtained from solving the eikonal equation with ray tracing and related techniques.
\end{remark}

%!TEX root = draft.tex

\section{Learning Local Dominant Ray Directions} \label{section:NMLA}

In Section \ref{section:ray_FEM} we use geometric optics to provide the motivation for the ray-FEM by building an adaptive approximation space that incorporates ray information specific to the underlying Helmholtz equation. However, the ray directions, which depend on the medium and source distribution, are unknown quantities themselves, hence they need to be computed or estimated. One way is to compute the global phase function, by either ray tracing or solving the eikonal equation, and take its gradient. As discussed in the introduction, computing the global phase function in a general varying medium can be difficult.

In the present paper we propose a totally different approach. This novel approach is based on learning the dominant ray directions by probing the same medium with the same source using a relative low-frequency wave. To be more specific, we first solve the Helmholtz equation \eqref{eq:Helmholtz} with the same speed function $c(\x)$, right hand side $f(\x)$ and boundary conditions but with a relative low-frequency  $\wt{\omega} \sim \sqrt{\omega}$ on a mesh with size $h=\cO(\wt{\omega}^{-2}) =\cO(\omega^{-1})$ with a standard finite element method, which is quasi-optimal in that regime. Then the local dominant ray directions are estimated based on the computed low-frequency wave field. The key point is that the low-frequency wave has probed the medium specific to the problem globally and only local dominant ray directions need to be learned, which allows us to handle multiple arrivals of wave fronts seamlessly. In particular, we use numerical micro-local analysis (NMLA), which is simple and robust, in this work to extract the dominant ray directions locally. However, this is a signal processing task that can be accomplished using other methods such as Prony's method \cite{Carriere_Moses:High_resolution_radar_target_modeling_using_a_modified_Prony_estimator}, Pisarenko's method \cite{Pisarenko_method}, MUSIC \cite{Schmidt:MUSIC}, matrix pencil \cite{Hua_Sarkar:Matrix_pencil_method_for_estimating_parameters_of_exponentially_damped/undamped_sinusoids_in_noise}, among many others.

\subsection{NMLA} \label{subsection:NMLA}
% !!!! I rewrote the paragraph I still don't like it...
In this subsection, for the sake of completness, we provide a brief introduction to NMLA developed in \cite{Benamou:NMLA_harmonic_wavefields, Benamou:NMLA_revisited}. If we suppose that a wave field is a weighted superposition of plane waves with the same wave number, but propagating in different directions, then the aim of NMLA is to extract from samples of the wave field, the directions and the weights. In the sequel we use a 2D example to illustrate the method.

%The aim of NMLA is to extract the local plane wave components, both their amplitudes and directions, from the measurement of a wave field that is a superposition of several plane waves with the same frequency (!!!!!This seems to be odd). We use a 2D example to illustrate the method.

Suppose that a wave field, denoted by $u(\x)$, is composed of $N$ plane waves around an observation point $\x_0$,
\begin{equation}
u(\x)= \sum_{n=1}^{N} B_{n}e^{ik(\x-\x_0)\cdot\widehat{\textbf{d}}_n}, \quad | \widehat{\textbf{d}}_n |=1.
\label{eq:plane-wave}
\end{equation}
We suppose that we can sample the wave field, $u(\x)$, and its derivative on a circle $S_r(\x_0)$ centered at $\x_0$ with radius $r$. The wave field can be written under the model assumption in \eqref{eq:plane-wave} as
\begin{equation} \label{eq:model_NMLA}
u(\textbf{x}_0 + r\widehat{\textbf{s}})= \sum_{n=1}^{N} B_{n}e^{i\alpha \widehat{\textbf{s}}\cdot\widehat{\textbf{d}}_n}, \quad \alpha=kr, ~\widehat{\textbf{s}}\in \mathbb{S}^1.
\end{equation}
Furthermore, define implicitly the angle variable $\theta = \theta(\widehat{\textbf{s}})$; accordingly we denote $\theta_n = \theta(\widehat{\textbf{d}}_n)$ the angles to recover, and $\x(\theta) = \textbf{x}_0 + r \widehat{\textbf{s}}(\theta)$. Using the angle based notation we sample the impedance quantity
\begin{equation}
U(\theta) =   \frac{1}{ik} \partial_{r} u(\textbf{x}(\theta)) + u(\textbf{x}(\theta)),
\end{equation}
which removes any possible ambiguity due to resonance, and improves the robustness to noise for solutions to the Helmholtz equation \cite{Benamou:NMLA_revisited}, on the circle $S_r(\x_0)$. Then we apply the filtering operator $\mathcal{B}$ to the impedance quantity
\begin{equation} \label{eq:filtered_data_1}
\mathcal{B}U(\theta) := \frac{1}{2L_{\alpha}+1}  \sum \limits_{l=-L_{\alpha}}^{L_{\alpha}}\frac{(\mathcal{F}U)_l e^{il\theta}}{(-i)^l (J_l(\alpha)-i J_l'(\alpha))} ,
\end{equation}
where $L_{\alpha} = \max (1, [\alpha], [\alpha + (\alpha)^{\frac{1}{3}} -2.5])$, $J_l$ is the Bessel functions of order $l$, $J'_l$ is their derivatives and
\begin{equation}
    \left( \mathcal{F}U\right)_l = \frac{1}{2\pi} \int_0^{2\pi} U(\theta) e^{-il \theta} d \theta
\end{equation}
is the $l$-th Fourier coefficient of $U$.
It is shown in \cite{Benamou:NMLA_revisited} that
\begin{equation} \label{eq:filtered_data_2}
\mathcal{B}U(\theta) = \sum_{n=1}^{N}B_nS_{L_{\alpha} }(\theta - \theta_n),
\end{equation}
where $ S_{L}(\theta) = \frac{\sin([2L+1]\theta/2)}{[2L+1]\sin(\theta/2)} $.
As a consequence, we have that if $\alpha=kr\rightarrow \infty$ then
\begin{equation}
\lim_{\alpha \rightarrow \infty}\mathcal{B} U(\theta) = \left\{ \begin{array}{ll}
 B_n, & \quad \textrm{if $\theta = \theta_n$ (or $\widehat{\textbf{s}} = \widehat{\textbf{d}}_n $ );}\\
 0, & \quad \textrm{otherwise}.
  \end{array} \right.
\end{equation}
Then it is possible to obtain the directions and the amplitudes by picking the peaks in the filtered data in \eqref{eq:filtered_data_2}.

However, for applications, the measured data is never a perfect superposition of plane waves; therefore, we provide, for completeness,  stability and error estimates for NMLA from \cite{Benamou:NMLA_revisited} in \ref{appendix-a}. In principle, as long as the perturbation is relative small with respect to the true signal, the estimation error is $\cO(\frac{1}{kr})$. In other words, the larger the radius of the circle compared to wavelength the more accurate the estimation is.

In our application, our data is the numerical solution of the Helmholtz equation. In addition to noises and numerical errors, there are perturbations due to two model errors :
\begin{itemize}
\item the geometric optics ansatz has an asymptotic error of order $\cO(\omega^{-1})$ (see \eqref{eq:plane_wave_expansion});
\item in the geometric ansatz the wave field at a point is a superposition of curved wave fronts. In particular, the curvature of the wave fronts results in a compromise in the choice of the radius of the sampling circle to be of order $\cO(\omega^{-\frac{1}{2}})$ for the NMLA in order to achieve the stability and the minimal error of order $\cO(\omega^{-\frac{1}{2}})$.

\end{itemize}
Detailed analysis is provided in \ref{appendix-b}. Below is a summary of the NMLA algorithm.

\alglanguage{pseudocode}
\begin{algorithm}  [H]
\caption{NMLA}
\label{alg:NMLA}
\begin{algorithmic}[1]
\Function {$ \textbf{d}_{\omega} $ = NMLA} {$\x_0, \mathcal{T}_h, \omega, h, c, \textbf{u} $}
\State choose $r \sim \omega^{-\frac{1}{2}}$   \Comment{Radius for the sampling circle}
\State choose $M \sim \omega r$  \Comment{Number of sampling points}
\State $\Delta \theta  = 2\pi/ M$, \Comment{Angular discretization}
\For {$ \theta = 0: \Delta \theta: 2\pi$}
\State $\x(\theta) = \x_0 + r\widehat{\textbf{s}}(\theta)$
\State $U(\theta) =   \frac{\omega}{ic(\x_0)} \partial_{r} u(\textbf{x}(\theta)) + u(\textbf{x}(\theta)) $ \Comment{Sample impedance data}
\EndFor
\State $F(\theta) = \mathcal{B}U(\theta)$  \Comment{Apply the filter \eqref{eq:filtered_data_1} }
\State $\boldsymbol{\theta}_{est} = \emph{SharpPeakLocations}(\theta, F(\theta))$
\State $\textbf{d}_{\omega} = \textbf{d}(\boldsymbol{\theta}_{est} )$
\EndFunction
\end{algorithmic}
\end{algorithm}

\subsection{Approximation property of numerical ray-FEM}
\label{subsec:numerical-ray-FEM}
% shall we make it shorter?

In this section we incorporate the errors in the estimation of the ray directions into the approximation error for the ray-FEM method, in which ray directions are first estimated by using Algorithm \ref{alg:NMLA} to the solution of the Helmholtz equation with a relative low-frequency. The estimated ray directions are then used to generate the approximation space.
With the same assumptions as in Section \ref{section:approximation_L2_norm}, we estimate an upper bound on
\begin{equation}
 \inf_{u_h \in V_{Ray}^h(\mathcal{T}_h)} \Vert u - u_h \Vert_{L^2(\Omega)},
\end{equation}
 when the ray-FEM space, $V_{Ray}^h(\mathcal{T}_h)$, is constructed using the estimated ray directions from high-frequency waves by NMLA.

From \ref{appendix-b}, the error estimation of dominant ray directions is $\mathcal{O}(\omega^{-1/2})$. The numerical ray-FEM space $V_{Ray}^h(\mathcal{T}_h)$ is defined similarly as $V_{Ray}(\mathcal{T}_h)$ with the exact ray directions $\lbrace \widehat{\textbf{d}}_j \rbrace $ replaced by the ones $\lbrace \widehat{\textbf{d}}_j^h \rbrace $ estimated by NMLA and $\vert \widehat{\textbf{d}}_j - \widehat{\textbf{d}}_j^h \vert \sim \mathcal{O}(\omega^{-1/2})$.

We denote by
\begin{equation} \label{nodal_interpolation_h}
u_I^h = \sum_{j=1 }^{N_h} A(\textbf{x}_j)\varphi_j(\textbf{x}) e^{i\omega [ \phi(\textbf{x}_j) +  1/c(\textbf{x}_j) \widehat{\textbf{d}}_j^h  \cdot ( \textbf{x} - \textbf{x}_j) ]}
\end{equation}
the nodal interpolation of the solution in $V_{Ray}^h(\mathcal{T}_h)$. Then we have
\begin{displaymath}
\begin{array}{ll}
\Vert u_I - u_I^h \Vert_{L^2(\Omega)}
& = \Vert \sum_{j=1 }^{N_h} A(\textbf{x}_j)\varphi_j(\textbf{x})
e^{i\omega  \phi(\textbf{x}_j)} ( e^{i\omega \nabla \phi(\textbf{x}_j) \cdot ( \textbf{x} - \textbf{x}_j) } - e^{i\omega  /c(\textbf{x}_j) \widehat{\textbf{d}}_j^h  \cdot ( \textbf{x} - \textbf{x}_j) } )
\Vert_{L^2(\Omega)} \\[1.0ex]

& \leq  \sum_{j=1 }^{N_h} \Vert A  \Vert_{L^{\infty}(\Omega)} \Vert e^{i\omega  /c(\textbf{x}_j) \widehat{\textbf{d}}_j  \cdot ( \textbf{x} - \textbf{x}_j)} - e^{i\omega /c(\textbf{x}_j) \widehat{\textbf{d}}_j^h  \cdot ( \textbf{x} - \textbf{x}_j) } \Vert_{L^2(S_j)}
\\[1.5ex]
& \lesssim \sum_{j=1 }^{N_h} \Vert A  \Vert_{L^{\infty}(\Omega)} \omega h \Vert c^{-1} \Vert_{L^{\infty}(\Omega)} \vert \widehat{\textbf{d}}_j - \widehat{\textbf{d}}_j^h \vert  |S_j|  \\[1.5ex]
& \lesssim  \omega^{1/2} h \Vert A  \Vert_{L^{\infty}(\Omega)}  \Vert c^{-1}  \Vert_{L^{\infty}(\Omega)}. \\
\end{array}
\end{displaymath}

Hence,
\begin{equation}
\begin{array}{ll}
\inf_{u_h \in V_{Ray}^h(\mathcal{T}_h)} \Vert u - u_h \Vert_{L^2(\Omega)} &\leq \Vert u - u_{I}^h \Vert_{L^2(\Omega)} \leq \Vert u - u_{I} \Vert_{L^2(\Omega)} + \Vert u_I - u_{I}^h \Vert_{L^2(\Omega)}\\ [1.5ex]
 & \lesssim h^2 \vert A \vert_{H^2(\Omega)} +  \omega h^2 \Vert A  \Vert_{L^{\infty}(\Omega)}  \Vert \nabla^2 \phi  \Vert_{L^{\infty}(\Omega)} \\ [1.5ex]
 & \quad   +  \omega^{1/2} h \Vert A  \Vert_{L^{\infty}(\Omega)}  \Vert c^{-1}  \Vert_{L^{\infty}(\Omega)} +  O(\omega^{-1}).
\end{array}
\end{equation}
Or, more compactly, we have that
\begin{equation}  \label{eq:approximation_property_numerical_ray}
    \inf_{u_h \in V_{Ray}^h(\mathcal{T}_h)} \Vert u - u_h \Vert_{L^2(\Omega)}  = \cO(h^2 + \omega h^2 + \omega^{1/2}h + \omega^{-1} ).
\end{equation}

Comparing with \eqref{eq:approximation_property_exact_ray} and \eqref{eq:approximation_property_numerical_ray}, the error in the estimation of dominant ray directions due to NMLA leads to the extra term $\omega^{1/2}h$, which is the leading order in the high-frequency regime. Specifically, if   $\omega h = \cO(1)$, then we have
\begin{equation}  \label{eq:approximation_property_numerical_ray_final}
    \inf_{u_h \in V_{Ray}^h(\mathcal{T}_h)} \Vert u - u_h \Vert_{L^2(\Omega)}  = \cO(\omega^{-1/2}).
\end{equation}

% texfile with the algorithms
%!TEX root = draft.tex

\section{Algorithms} \label{Section:Algorithms}

In this section we provide the full algorithm for ray-FEM including a fast iterative solver based on a modification of the method of polarized traces for the resulting linear systems. In order to streamline the presentation and to make the algorithm easier to understand, we introduce several subroutines, which the main algorithm relies on.

We can separate the full algorithm into three conceptual stages:
\begin{enumerate}
\item probing the medium by solving a relative low-frequency Helmholtz equation with the standard FEM,
\item learning the dominant ray directions from low-frequency probing wave field by NMLA,
\item solving the high-frequency Helmholtz equation in ray-FEM space.
\end{enumerate}

If necessary the second stage can be applied to the high-frequency wave field computed in stage 3 to improve the estimation of dominant ray directions and then repeat stage 3.

We remind the reader that the ultimate objective of the algorithm presented in this paper (i.e., Algorithm \ref{alg:Iter_Ray_FEM}) is to solve the Helmholtz equation \eqref{eq:Helmholtz} at frequency $\omega$ with a total $\cO(\omega^d)$ (up to poly-logarithmic factors) computational complexity. In order to achieve this objective, we discretize the PDE with a mesh size $h = \cO(\omega^{-1})$, which leads to a total of $\cO(\omega^d)$ number of degrees of freedom and a sparse linear system with  $\cO(\omega^d)$ number of nonzeros. Then we develop a fast iterative solver with quasi-linear complexity to solve the resulting linear system after discretization. Below is a more detailed description of the three stages. Finally, following the notation defined in the prequel, we denote the triangular mesh by $\mathcal{T}_h$.

\subsection{Probing}
We first solve the Helmholtz equation with frequency $\widetilde{\omega} \sim \sqrt{\omega}$ in the same medium and with the same source on $\mathcal{T}_h$. The low-frequency problem is solved using the standard finite element method (S-FEM) with linear elements as prescribed by Algorithm \ref{alg:S-FEM}.
%\alglanguage{pseudocode}
\begin{algorithm}[H]
\caption{Standard FEM Helmholtz Solver}
\label{alg:S-FEM}
\begin{algorithmic}
\Function {$\textbf{u}_{\omega, h}$ = S-FEM} {$\omega, h, c, f, g$}
\For {$i,j = 1:N_h$}
 \State $\textbf{H}_{i,j} = \mathcal{B}(\varphi_i, \varphi_j)$   \Comment{Assemble Helmholtz matrix}
 \State $\textbf{b}_j = \mathcal{F}(\varphi_j)$ \Comment{Assemble right-hand side}
\EndFor
 \State $\textbf{u}_{\omega, h} = \textbf{H}^{-1} \textbf{b} $ \Comment{Solve linear system}
\EndFunction
\end{algorithmic}
\end{algorithm}

Let $\textbf{u}_{\widetilde{\omega}, h}$ = S-FEM ($\widetilde{\omega}, h, c, f, g$) denote the S-FEM solution of the low-frequency Helmholtz equation on $\mathcal{T}_h$. Since $\tilde{\omega}^2 h = \mathcal{O}(1)$, S-FEM is quasi-optimal in the norm $\Vert \cdot \Vert_{ \mathcal{H}} := \Vert \nabla \cdot \Vert_{L^2} + k \Vert \cdot \Vert_{L^2} $ \cite{Melenk:Generalized_FEM}, and it has optimal $L^2$ error estimate \cite{Haijun_Wu_14:CIP_FEM_linear_version}.

\subsection{Learning}
\label{sec:learning}
Once the low-frequency problem has been solved, we extract the dominant ray directions from $\textbf{u}_{\widetilde{\omega}, h}$ using NMLA as described in Section \ref{subsection:NMLA} around each mesh node. We utilize the smoothness of the phase functions, and hence the smoothness of the ray directions field to reduce the computational cost. The reduction is achieved by restricting the learning of the dominant ray directions to vertices of a coarse mesh down-sampled from $\mathcal{T}_h$. Such coarse mesh is denoted by $\mathcal{T}_{h_c}$, where $h_c=\cO(\sqrt{h})$. The resulting dominant ray directions are then interpolated onto the fine mesh $\mathcal{T}_h$.

%(REWRITE paragraph :/)
Note that at each vertex of $\mathcal{T}_{h_c}$, the wave field $\textbf{u}_{\widetilde{\omega}, h}$ on the fine mesh $\mathcal{T}_h$ is used for NMLA to estimate the dominant ray directions. There are three sources of error in the learning stage we need to be aware of:

\begin{itemize}
\item numerical error of $\textbf{u}_{\widetilde{\omega}, h}$,
\item model error in geometric optics ansatz,
\item interpolation error.
\end{itemize}

The numerical error for $\textbf{u}_{\widetilde{\omega}, h}$ by Algorithm \ref{alg:S-FEM} in $L^2$ norm \cite{Haijun_Wu_14:CIP_FEM_linear_version} is $\cO(\wt{\omega}h^2 + \wt{\omega}^2h^2) = \cO(\omega^{-1})$,
which is negligible with respect to the model error in the geometric optics ansatz. The error introduced by the geometric optics approximation is $\cO(\widetilde{\omega}^{-\frac{1}{2}})$ as shown in Section \ref{subsection:NMLA} and \ref{appendix-b}. The error due to linear interpolation on $\mathcal{T}_{h_c}$ to get ray direction estimation at every vertex on $\mathcal{T}_h$ is $\cO(h_c^2)=\cO(h)=\cO(\omega^{-1})$, which is much smaller than the model error in geometric optics ansatz. Hence the overall error in ray direction estimation based on NMLA on $\textbf{u}_{\widetilde{\omega}, h}$  and interpolation is $\cO(\widetilde{\omega}^{-\frac{1}{2}})$. The dominant ray direction estimation algorithm is summarized in Algorithm \ref{alg:Ray_Learning}.
For each node $\x_j$ on mesh $\mathcal{T}_h$ the number of dominant ray directions is denoted by $n_j$, $\boldsymbol{d}_{\omega,h}^{j} =  \lbrace \textbf{d}_{\omega, h}^{j, l} \rbrace_{l = 1}^{n_{j}} $.

%\alglanguage{pseudocode}
\begin{algorithm}  [H]
\caption{Ray Learning }
\label{alg:Ray_Learning}
\begin{algorithmic}[1]
\Function {$ \lbrace \boldsymbol{ d}_{\omega,h}^j \rbrace_{j=1}^{N_h}$ = RayLearning} {$\omega, h, h_c, c, \textbf{u}_{\omega, h} $}
\For  {$j = 1: N_{h_c}$ }
\State $\boldsymbol{d}_{\omega,h_c}^{j} = \emph{NMLA} ( \x_j^c, \mathcal{T}_h, \omega, h, c, \textbf{u}_{\omega,h})$
\EndFor
\State $\lbrace \boldsymbol{ d}_{\omega,h}^j \rbrace_{j=1}^{N_h} = \emph{Interpolation} (\mathcal{T}_{h_c}, \mathcal{T}_h, \lbrace \boldsymbol{d}_{\omega, h_c}^{j} \rbrace_{j = 1}^{N_{h_c}} )$
\EndFunction
\end{algorithmic}
\end{algorithm}

% this is what Jun used to have, I haven't erased.

%In the present algorithm the probing will be performed in both high-frequency wavefield $\textbf{u}_{\omega, h}$ and low-frequency wavefield $\textbf{u}_{\widetilde{\omega}, h}$ on the same fine and coarse mesh. We denote by $ \lbrace \boldsymbol{ d}_{\omega,h}^j \rbrace_{j=1}^{N_h}$ and $ \lbrace \boldsymbol{ d}_{\widetilde{\omega},h}^j \rbrace_{j=1}^{N_h}$ the resulting ray directions, i.e., $ \lbrace \boldsymbol{ d}_{\omega,h}^j \rbrace_{j=1}^{N_h} = \emph{RayProbing}  \lbrace \omega, h, h_c, \textbf{u}_{\omega, h} \rbrace$, $ \lbrace \boldsymbol{ d}_{\widetilde{\omega},h}^j \rbrace_{j=1}^{N_h} =  \emph{RayProbing} \lbrace \widetilde{\omega}, h, h_c, \textbf{u}_{\widetilde{\omega}, h} \rbrace $.

%\begin{remark} {\label{remark:rayprobing}}
%Since NMLA requires that the observation point is far from the source, for mesh nodes near the source, we simply use the analytical propagating direction in homogeneous medium to approximate the true ray direction.
%Also NMLA needs the wavefield data $\textbf{u}_{\omega, h} $ on the observation circles centered at all mesh nodes, the domain of $\textbf{u}_{\omega, h} $ should be slightly larger than the domain of ray probing field.
%The choice of coarse mesh size $h_c$ and interpolation method will be described in section \ref{subsection:complexity_ray_probing}.
%\end{remark}

\subsection{High-frequency solver}

Once the dominant ray directions on $\mathcal{T}_h$ have been computed, we can construct the ray-FEM space $V_{Ray}(\mathcal{T}_h)$ and solve the high-frequency Helmholtz equations following \eqref{R-FEM}, which is implemented in Algorithm \ref{alg:Ray_FEM}.

%\alglanguage{pseudocode}
\begin{algorithm}
\caption{Ray-FEM Helmholtz Solver }
\label{alg:Ray_FEM}
\begin{algorithmic}[1]
\Function {$\textbf{u}_{\textbf{d}_{\omega}, h}$ = Ray-FEM} {$\omega, h, c, f, g, \lbrace \boldsymbol{ d}_{\omega,h}^j \rbrace_{j=1}^{N_h} $}
\State $N_{dof} = 0$
\For {$j = 1: N_h, l = 1:n_j$}
\State $N_{dof} = N_{dof} + 1$, $m = N_{dof}$
\State $\psi_m(\textbf{x}) = \varphi_j (\x) e^{iw/c(\x_j) \textbf{d}_{j,l}\cdot \textbf{x}}$
\Comment{Construct ray-FEM basis functions}
\State $\widehat{\psi} = \psi_m(\textbf{x}_j) $
\Comment{Nodal values of ray-FEM basis functions}
\EndFor
\For {$m,n = 1:N_{dof}$}
 \State $\textbf{H}_{m,n} = \mathcal{B}(\psi_m, \psi_n)$   \Comment{Assemble Helmholtz matrix}
 \State $\textbf{b}_n = \mathcal{F}(\psi_n)$ \Comment{Assemble right-hand side}
\EndFor
\State $\textbf{v} = \textbf{H}^{-1} \textbf{b} $ \Comment{Coefficients of ray-FEM basis functions}
\State $\textbf{u}_{\textbf{d}_{\omega}, h} = \textbf{v} \cdot \boldsymbol{ \widehat{  \psi }} $ \Comment{Ray-FEM solution on mesh nodes}
\EndFunction
\end{algorithmic}
\end{algorithm}

\iffalse
The input ray directions $\lbrace \boldsymbol{ d}_{\omega,h}^j \rbrace_{j=1}^{N_h}$ in above algorithm is not only restricted to ray probing for high-frequency waves, they can also be replaced by low-frequency ray estimation $\lbrace \boldsymbol{d}_{\widetilde{\omega},h}^j \rbrace_{j=1}^{N_h}$ or exact ray directions $\lbrace \boldsymbol{ d}_{ex}^j \rbrace_{j=1}^{N_h}$, the corresponding ray-FEM solutions are denoted by $\textbf{u}_{\textbf{d}_{\widetilde{\omega}}, h}$ or $\textbf{u}_{\textbf{d}_{ex}, h}$, respectively.
\fi

In general, the more accurate of the dominant ray directions learned by NMLA, the more accurate ray-FEM high-frequency solution can be computed by Algorithm \ref{alg:Ray_FEM}. From section \ref{sec:learning}, the accuracy order of learning stage from low-frequency wave field is $\cO(\wt{\omega}^{-\frac{1}{2}})$, and following the error analysis of section \ref{subsec:numerical-ray-FEM}, the consequent ray-FEM solution has the same order of accuracy. However, we can apply the learning stage to the computed high-frequency wave field to improve approximation for dominant ray directions. Moreover, the improved ray information can be used to improve the approximation for high-frequency wave field computed by ray-FEM. This process can be continued iteratively to further improve the approximations of the ray directions and the high-frequency wave field. Therefore, we develop an iterative ray-FEM Helmholtz solver for high-frequency $\omega$ in Algorithm \ref{alg:Iter_Ray_FEM}.

% previous version:
%In general, an iterative Ray-FEM Helmholtz solver for high-frequency $\omega$ can be designed as follows. We start with a relative low-frequency wave $\widetilde{\omega} \sim \omega^{1/2}$  with the same source distribution to probe and learn the medium and then extract dominant ray directions. These learned ray information will be used to compute the Ray-FEM high-frequency solution on a mesh with size $\omega h = \mathcal{O}(1)$ due to improved stability and accuracy compared to using standard Galerkin basis. Furthermore, we iteratively use NMLA and Ray-FEM with newly computed high-frequency Helmholtz solution and ray information to improve approximations for both ray directions and high-frequency wave field.

%\alglanguage{pseudocode}
\begin{algorithm}
\caption{Iterative Ray-FEM High-Frequency Helmholtz Solver}
\label{alg:Iter_Ray_FEM}
\begin{algorithmic}[1]
\Function {$\textbf{u}_{\textbf{d}_{\omega}, h}$ = IterRay-FEM} {$\omega, c, f, g$}
\State $\widetilde{\omega} \sim \sqrt{\omega}$, $h \sim \omega^{-1}$, $h_c \sim \omega^{-\frac{1}{2}} $
\State $\textbf{u}_{\widetilde{\omega}, h} = \emph{S-FEM} (\widetilde{\omega}, h, c, f, g)$ \Comment{Low-frequency waves}
\State $\lbrace \boldsymbol{d}_{\widetilde{\omega},h} \rbrace = \emph{RayLearning} (\widetilde{\omega}, h, h_c, c,  \textbf{u}_{\widetilde{\omega}, h} )$
\Comment{Low-frequency ray learning }
\State $\textbf{u}_{\textbf{d}_{\widetilde{\omega}}, h} = \emph{Ray-FEM} (\omega, h, c, f, g, \lbrace \boldsymbol{d}_{\widetilde{\omega},h} \rbrace)$ \Comment{High-frequency waves}
\State  $tol = 1$, $niter = 0$, $\textbf{u}_{\omega,h}^1 = \textbf{u}_{\textbf{d}_{\widetilde{\omega}}, h}$
\While {$tol > \epsilon $ or $niter> max\_iter$}
         \State $\lbrace  \boldsymbol{\textbf{d}}_{\omega,h} \rbrace = \emph{RayLearning} (\omega, h, h_c, c, \textbf{u}_{\omega,h}^1 )$  \Comment{High-frequency ray learning }
         \State $\textbf{u}_{\textbf{d}_{{\omega}}, h} = \emph{Ray-FEM} (\omega, h, c, f, g, \lbrace \boldsymbol{\textbf{d}}_{\omega,h} \rbrace)$, $\textbf{u}_{\omega,h}^2 = \textbf{u}_{\textbf{d}_{{\omega}}, h}$
         \State $tol = \Vert \textbf{u}_{\omega,h}^1 - \textbf{u}_{\omega,h}^2 \Vert_{L^2(\Omega)} / \Vert \textbf{u}_{\omega,h}^2 \Vert_{L^2(\Omega)}$
         \State  $niter = niter+1$, $\textbf{u}_{\omega,h}^1 = \textbf{u}_{\omega,h}^2$
      \EndWhile
% \State $\textbf{u}_{\textbf{d}_{\omega}, h} = \textbf{u}_{\omega,h}^2$
\EndFunction
\end{algorithmic}
\end{algorithm}

\begin{remark} {\label{remark:rayprobing}}
Extensive numerical experiments and \ref{appendix-a} suggest that the NMLA process in learning dominant ray directions is remarkably stable even for noisy plane wave data. Hence, the iterative process in Algorithm \ref{alg:Iter_Ray_FEM} usually needs very few iterations to reach the desired accuracy. Typically, we only need $1$ or $2$ iterations in our numerical tests.
\end{remark}

\subsection{Fast linear solver} \label{subsection:fast_linear_solver}

To achieve the overall complexity mentioned in the introduction, it is necessary to solve the linear system resulting from both standard finite element methods and ray-FEM, which we write in a generic form as
\begin{equation} \label{eq:discrete_Helmholtz_equation}
\H\u = \f,
\end{equation}
in linear complexity (up to poly-logarithmic factors). This solver is, in fact, the bottleneck of Algorithms \ref{alg:S-FEM} and \ref{alg:Ray_FEM}.

For a smooth medium, this can be achieved by modifying the method of polarized traces \cite{ZepedaDemanet:the_method_of_polarized_traces}, of which we provide a brief review here. For further details we refer the interested readers to \cite{ZepedaDemanet:the_method_of_polarized_traces}. The method of polarized traces is a domain decomposition method that encompasses the following aspects
\begin{itemize}
\item  layered domain decomposition;
\item  absorbing boundary conditions between subdomains implemented via PML \cite{Berenger:PML};
\item  transmission conditions issued from a discrete Green's representation formula;
\item  efficient preconditioner arising from localization of the waves via an incomplete Green's formula.
\end{itemize}

The first two aspects can be effortlessly implemented. Consider a layered partition of $\Omega$ into $L$ slabs, or layers $\{ \Omega^{\ell} \}_{\ell =1}^{L}$. Define $f^{\ell}$ as the restriction of $f$ to $\Omega^{\ell}$, i.e., $f^{\ell} = f\chi_{\Omega^{\ell}}$; and define the local Helmholtz operators as
\begin{equation} \label{eq:local_Helmholtz}
\cH^{\ell} u  :=  \left( -\triangle  - m \omega^2  \right)  u   \qquad  \text{in } \Omega^{\ell},
\end{equation}
with absorbing boundary conditions implemented via PML between slabs.

The method of polarized traces aims to solve the global linear system in \eqref{eq:discrete_Helmholtz_equation} by solving the local systems $\H^{\ell}$, which are the discrete version of \eqref{eq:local_Helmholtz}.

In order to solve the global system, or in this case to find a good approximate solution, we need to ``glue'' the subdomains together, which is achieved via a discrete Green's integral formula deduced by imposing discontinuous solutions.

In the original formulation of the method of polarized traces \cite{ZepedaDemanet:the_method_of_polarized_traces}, the Green's representation formula was used to build a global surface integral equation (SIE) at the interfaces between slabs. The SIE was solved using an efficient preconditioner coupled with a multilevel compression of the discrete kernels to accelerate the on-line stage of the algorithm. The original algorithm had superlinear off-line complexity which was amortized among a large number of right-hand sides, which represents a typical situation in explorations geophysics.

In the context of the present paper, the linear systems issued from the ray-based FEM depends on the source distribution, making it impossible to amortize a super-linear off-line cost. In order to reduce the off-line cost we use a matrix-free formulation (see Chapter 2 in \cite{Zepeda_Nunez:Fast_and_scalable_solvers_for_the_Helmholtz_equation}) with a domain decomposition in thin layers. In this case, the cost per iteration is linear on the number of degrees of freedom, depending on the grow or the auxiliary degrees of freedom corresponding to the PML's. Finally the convergence is normally achieved in $\mathcal{O}(\log{\omega})$ iterations, as it will be shown in the sequel.

% texfile with the complexity
%!TEX root = draft.tex

\section{Complexity} \label{section:complexity}

In this section we provide an overall computational complexity count of our algorithm for the high-frequency Helmholtz equation \eqref{eq:Helmholtz} in terms of $\omega$. Our count includes learning ray directions by NMLA and the linear solver for the discretized systems from both standard FEM and ray FEM for low-frequency and high-frequency Helmholtz equation respectively.
The summary of the complexity of the method is given in Table \ref{tab:complexity_all_methods}.

\subsection{Ray direction learning} \label{subsection:complexity_ray_probing}

As described in Section \ref{sec:learning},  Algorithm \ref{alg:Ray_Learning} applies NMLA to computed wave fields with low-frequency $\widetilde{\omega} \sim \sqrt{\omega}$ or high-frequency $\omega$ and first estimate ray directions at vertices on a downsampled coarse mesh $\mathcal{T}_{h_c}$ and then interpolate the ray directions to the vertices on a fine mesh $\mathcal{T}_{h}$. We remind the reader the following scalings: $h=\cO(\omega^{-1})$, $h_c=\cO(\sqrt{h})=\cO(\omega^{-\frac{1}{2}})$ . These scalings allow us to strike a balance among the number of observation points at which NMLA is used to estimate ray directions, the radius of sampling circle, and the corresponding number of sampling points on the circle to resolve the wave field to reach the optimal accuracy of NMLA with desired total computational complexity.

We estimate the dominant computational cost which is applying NMLA to the high-frequency wave field to get ray direction estimation on $\mathcal{T}_{h}$ in 2D as an illustration of Table \ref{tab:complexity_ray_probing}. It is shown in \ref{appendix-b}, the least error that can be achieved by NMLA is $\cO(\omega^{-\frac{1}{2}})$ when radius $r$ of the sampling circle centered at an observation point is $\cO(\omega^{-\frac{1}{2}})$. Hence the number of points sampled on the circle  to resolve the wave field of frequency $\omega$ is $M_{\omega}=\cO(\omega^{\frac{1}{2}})$. Since NMLA is a linear filter based on a Fourier transform in angle space, the corresponding computational complexity is $\cO(M_{\omega}\log M_{\omega})$ \cite{Benamou:NMLA_harmonic_wavefields}. The number of observation points we need to perform NMLA is the number of vertices on the coarse mesh which is $\cO(h_c^{-2})=\cO(\omega)$. Hence the computation cost to obtain ray directions at all vertices on the coarse mesh by NMLA is $\cO(\omega^{\frac{3}{2}}\log \omega)$. Finally the ray directions estimated at the vertices on the coarse mesh by NMLA are interpolated to the fine mesh $\mathcal{T}_{h}$. Interpolation is a linear operation and hence its computation complexity is $\cO(\omega^2)$.

 Table \ref{tab:complexity_ray_probing} provides the complexity for each operation of NMLA, where $d$ is the dimension and $C_{NMLA}$, $C_{ray, h_c}$,  $C_{Int}$, and $C_{ray, h}$ is the computation complexity of NMLA at a single vertex, NMLA on the under sampled coarse mesh, interpolation of local ray directions to the fine mesh and full algorithm for learning local ray directions at frequency $\omega$ on mesh $\mathcal{T}_{h}$ respectively.

\begin{table}[h]
%title of the table
\centering % centering table
\begin{tabular}{||c ||c|c|c|c|c|c||} % creating eight columns
\hline\hline %inserting double-line
Frequency &  \quad $r$ \quad & $M_{\omega}$ & $C_{NMLA} $ & $C_{ray, h_c}$ & $C_{Int}$ & $C_{ray, h}$  \\
\hline
$\omega$   &  $\omega^{-\frac{1}{2}} $ & $\omega^{\frac{d-1}{2}} $ & $\omega^{\frac{d-1}{2}} \log \omega $  &  $\omega^{d-\frac{1}{2}} \log \omega $ & $\omega^d$ & $\omega^d$  \\
%\hline    % inserts single-line  [0.5ex]
%$\widetilde{\omega}$  &  $\widetilde{\omega}^{-\mu} $ & $\omega^{(1-\mu)(d-1)/2} $ & $\omega^{(1-\mu)(d-1)/2} \log \omega $   &  $\omega^{1 + (1-\mu)(d-1)/2} \log \omega$ & $\omega^d$ & $\omega^d$  \\
\hline \hline
\end{tabular}
\caption{Computational complexities of estimating ray directions on a mesh $\mathcal{T}_{h}$ with $h =\cO(\omega^{-1})$. }
\label{tab:complexity_ray_probing}
\end{table}

\subsection{Helmholtz Solver}

The most computationally intensive component in the whole ray-FEM algorithm is solving the linear systems after discretization of the Helmholtz equation. Algorithm \ref{alg:Iter_Ray_FEM} solves both $\textbf{u}_{\widetilde{\omega}, h} =  \emph{S-FEM} (\widetilde{\omega}, h, c, f, g)$ by S-FEM and $\textbf{u}_{\textbf{d}_{\omega}, h} = \emph{Ray-FEM}(\omega, h, c, f, g, \lbrace \boldsymbol{ d}_{\omega,h}^j \rbrace_{j=1}^{N_h} )$ by ray-FEM on the same mesh $\mathcal{T}_{h}$.  Each solver is composed of three steps: the assembling step, the setup step, and the iterative solve step.

Since the basis functions are locally supported, the resulting matrix is sparse. The complexity of the assembling step is of the same order as the degrees of freedom $N_h=\cO(\omega^d)$.

In the setup stage, the computational domain is decomposed into subdomains of thin layers whose width is comparable to the characteristic wavelength. The local problems in each subdomain are factorized using multifrontal method in $\cO(\sqrt{N_{h}})$ time (or  $\cO(\sqrt{N_{h}} \log^3 {N_{h}})$ time depending on the width of the auxiliary PML for each subdomain in term of the wavelength). Given that the layers are $\cO(1)$ elements thick, we have to factorize $\cO(\sqrt{N_{h}})$ subsystems, which results in a total  $\cO(N_{h})$ (or $\cO(N_{h} \log^3{N_{h}})$) asymptotic complexity for the setup step.

Finally, for the iterative solve step, each application of the preconditioner involves $6$ local solves per layer, each one performed with $\cO(\sqrt{N_{h}})$ ( or  $\cO(\sqrt{N_{h}} \log^2{N_{h}})$) complexity. Then given that we have $\cO(\sqrt{N_{h}})$ layers, then we have an overall $\cO(N_{h})$ (or $\cO(N_{h} \log{N_{h}})$) complexity per iteration. Extensive numerical experiments suggest that the number of iterations to converge is $\cO(\log{N_{h}})$ for the both high- and low-frequency solves for smooth media. Hence, the empirical overall complexity is $\cO(N_{h} \log{N_{h}})$ for the high-frequency solve and  $\cO(N_{h} \log^3{N_{h}})$ for the low-frequency one, which is summarized in the following table \ref{tab:complexity_all_methods}:

\begin{table}[h]
%title of the table
\centering % centering table
\begin{tabular}{||c ||c|c|c|c||} % creating eight columns
\hline\hline %inserting double-line
Methods & \quad S-FEM \quad  & \quad Learning \quad & \quad Ray-FEM \quad & Iterative Ray-FEM  \\
\hline
Frequency   &  $\sqrt{\omega} $  &  $\sqrt{\omega} \mbox{ or } \omega $ & $\omega $ & $\omega$\\
\hline    % inserts single-line  [0.5ex]
Complexity   & $\cO(\omega^d \log^3{\omega})$ &  $\cO(\omega^d)$ & $\cO(\omega^d \log{\omega})$ & $\cO(\omega^d \log^3{\omega})$  \\
\hline \hline
\end{tabular}
\caption{Overall computational complexities with mesh size $h = \cO(\omega^{-1})$.}
\label{tab:complexity_all_methods}
\end{table}

% little tex file with all the parts of the fast solver
%\input{fast_solver}

% tex file with all the numerical experiments
%!TEX root = draft.tex

\section{Numerical Experiments} \label{Section:Numerical_Experiments}

In this section we provide several numerical experiments to test the proposed ray-FEM and corroborate our claims. For all cases, the domain of interest is $\Omega = [-1/2,1/2]^2$ with different source terms and boundary conditions. $\Omega$ is discretized using a standard traingular mesh. The mass and stiffness matrices are assembled using a  high-order Gaussian quadrature rule to compute the integrals numerically\footnote{ Given the expression of the mass and stiffness matrices, which are polynomials times a plane wave, it is possible to compute the integral analytically \cite{Perugia_Pietra_Russo:PW_VEM}.}.

%We show three different kind of numerical experiments:
%\begin{itemize}
%\item source outside the domain: the experiments in Section \ref{section:source_outside} are meant show the convergence rate of the solution generated by the algorithm presented in this paper;
%\item source inside the domain: the experiments in Section \ref{section:source_inside} shows that the discretization is able to reduce the pollution error, and obtain better phases in the far field;
%\item using fast methods: in Section \ref{section:experiments_fast_methods} we present some examples in which fast methods are used to compute the solution to both low and high frequency problems.
%\end{itemize}

\subsection{Convergence tests} \label{subsection:convergence}
In the first test, the exact solution to the Helmholtz equation with Robin boundary condition is the wave field (normalized by the frequency $\omega$) corresponding to a point source outside the domain. It is given by
\begin{equation} \label{eq:exact_solution_1_point_source_outside}
  u_{ex}(x,y) = \sqrt{\omega}H_0^{(1)}(\omega\sqrt{(x-2)^2 + (y-2)^2}).
\end{equation}
Numerically we solve the Helmholtz equation \eqref{eq:Helmholtz} with $c(\x)\equiv 1$, source $f(\x)\equiv 0$ and exact impedance boundary data with mesh size such that the number of points per wave length (NPW) is 6 for different $\omega$ and test convergence for both the ray direction estimation by NMLA and the final numerical solution by ray-FEM.

%We compute the solution of the Helmholtz equation using Alg.~\ref{alg:Iter_Ray_FEM}
 First, a probing wave with low-frequency $\widetilde{\omega} = \sqrt{\omega}$ is solved by standard FEM. Then NMLA is applied to the low-frequency probing wave to get an estimation of the local dominant ray directions $\textbf{d}_{\widetilde{w}}$. Instead of using the regular NMLA for plane wave decomposition, we use NMLA with a curvature correction version \cite{Benamou:NMLA_revisited} to estimate the normal direction of a circular wave front. The local ray direction information is used in ray-FEM to produce the first numerical solution to the high frequency Helmholtz equation $u_{\textbf{d}_{\widetilde{w}}}$.

 We employ one more iteration in the framework of iterative ray-FEM by applying NMLA to $u_{\textbf{d}_{\widetilde{w}}}$ to get an improved local ray direction estimation $\textbf{d}_w$ and then use it again in ray-FEM to get a more accurate numerical solution to the high-frequency Helmholtz equation  $u_{\textbf{d}_w}$.

 Due to the curvature correction in NMLA, it can be shown that the error in ray direction estimation $\|\textbf{d}_{\widetilde{w}}-\textbf{d}_{ex}\|$ ($\|\textbf{d}_w-\textbf{d}_{ex}\|$) is $\cO (\widetilde{\omega}^{-1})$ ($\cO (\omega^{-1})$) by using an analysis similar to that in \cite{Benamou:NMLA_revisited} and \ref{appendix-b}. Using the estimate in Section \ref{subsec:numerical-ray-FEM}, one can show that the approximation error for numerical ray-FEM space is at least $\cO (\widetilde{\omega}^{-1})$ ($\cO (\omega^{-1})$) if $\textbf{d}_{\widetilde{w}}$ ($\textbf{d}_{w}$) is used.

\begin{figure}[h]
\begin{center}
{\includegraphics[trim = 0mm 0mm 0mm 0mm,clip, width=7cm ]{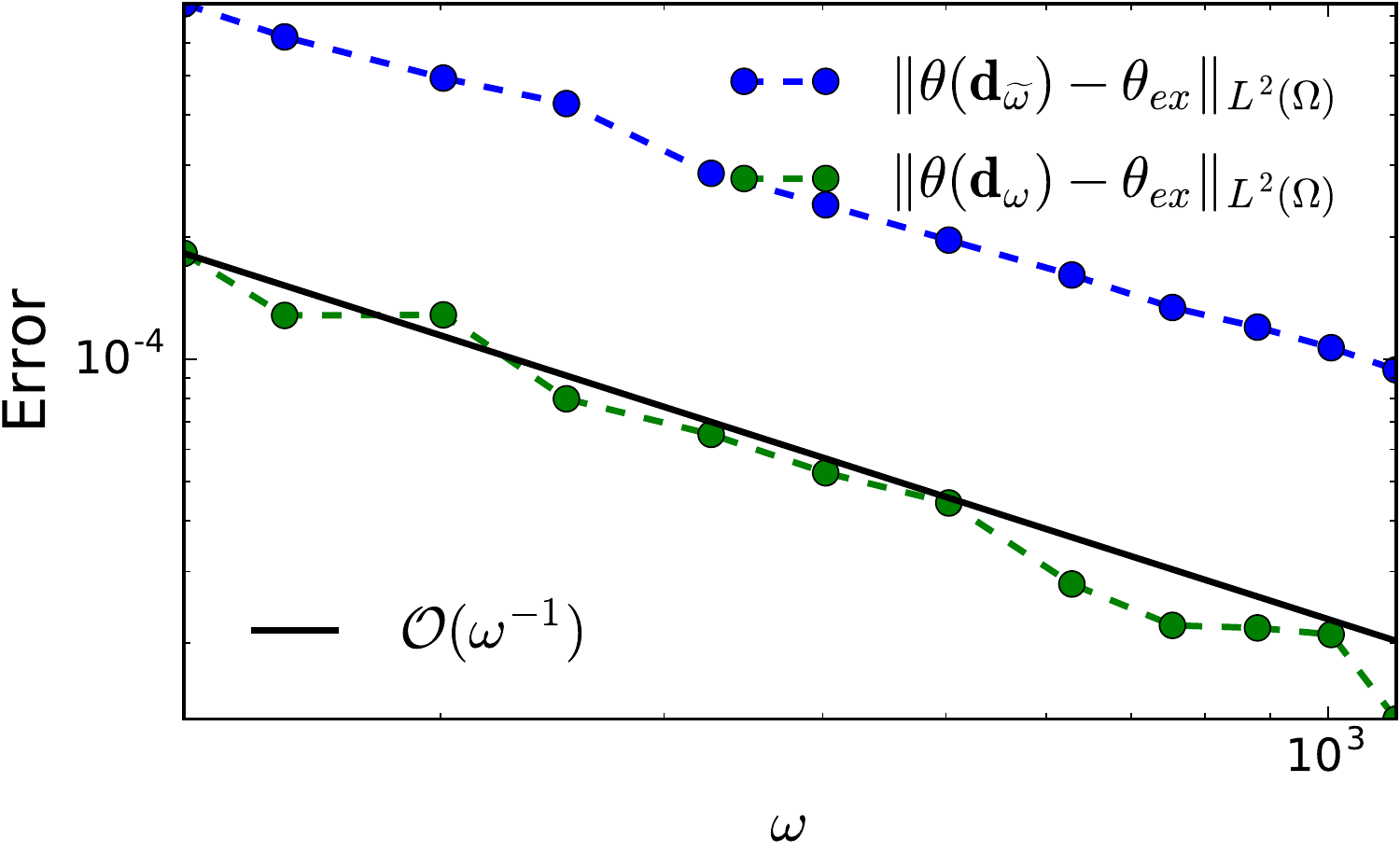}
\includegraphics[trim = 0mm 0mm 0mm 0mm,clip,width=7cm ]{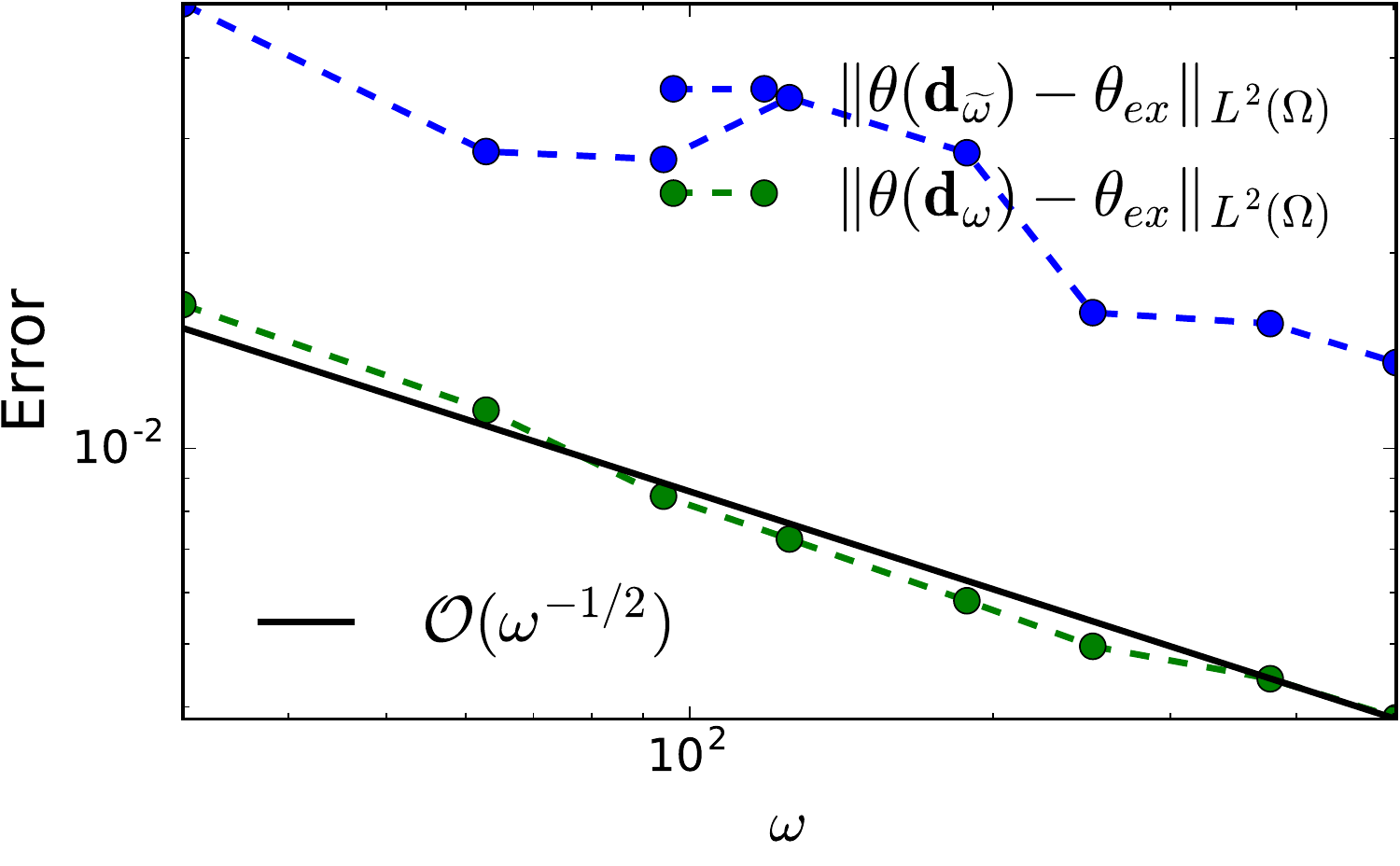}
\includegraphics[trim = 0mm 0mm 0mm 0mm,clip,width=7cm ]{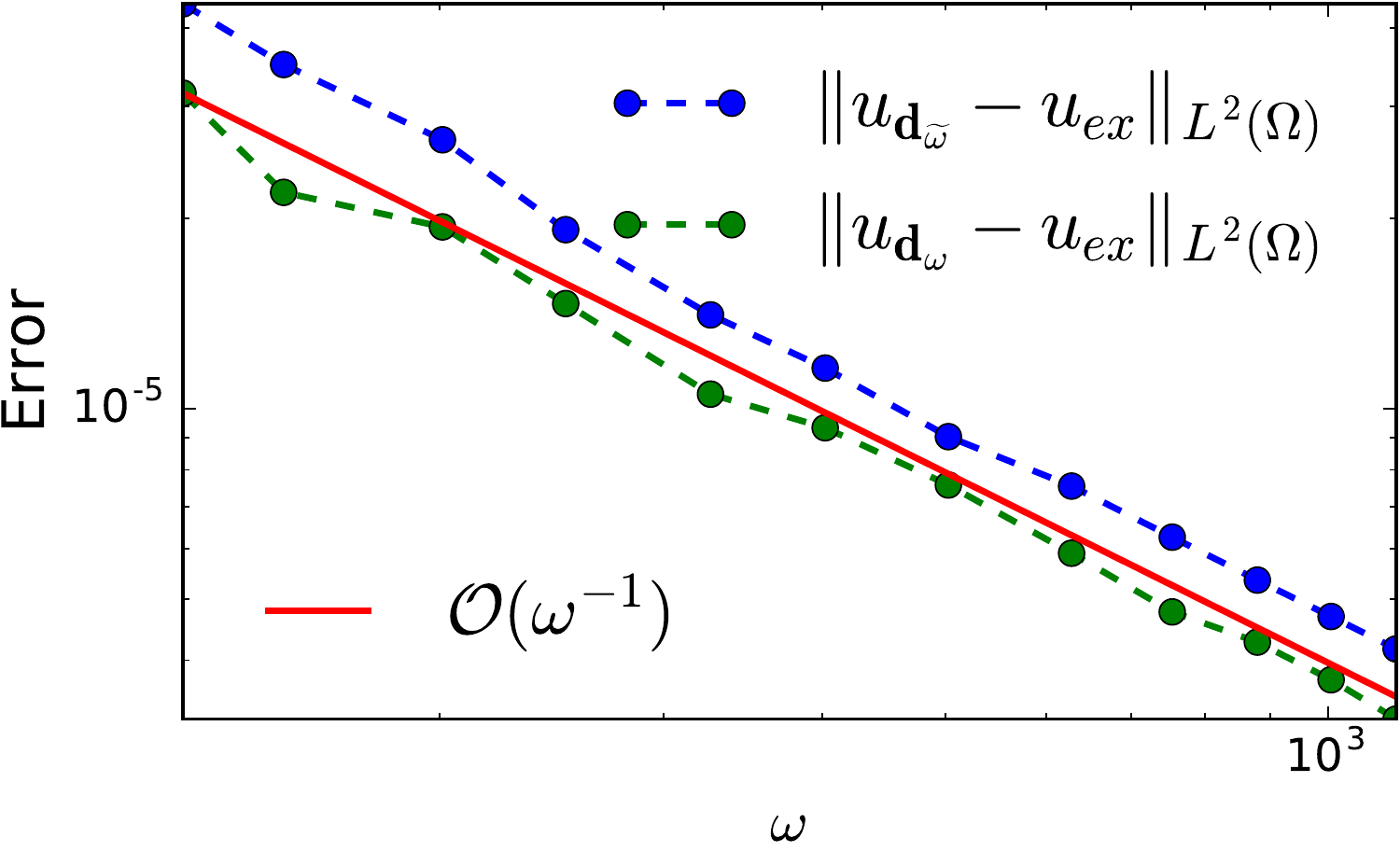}
\includegraphics[trim = 0mm 0mm 0mm 0mm,clip,width=7cm ]{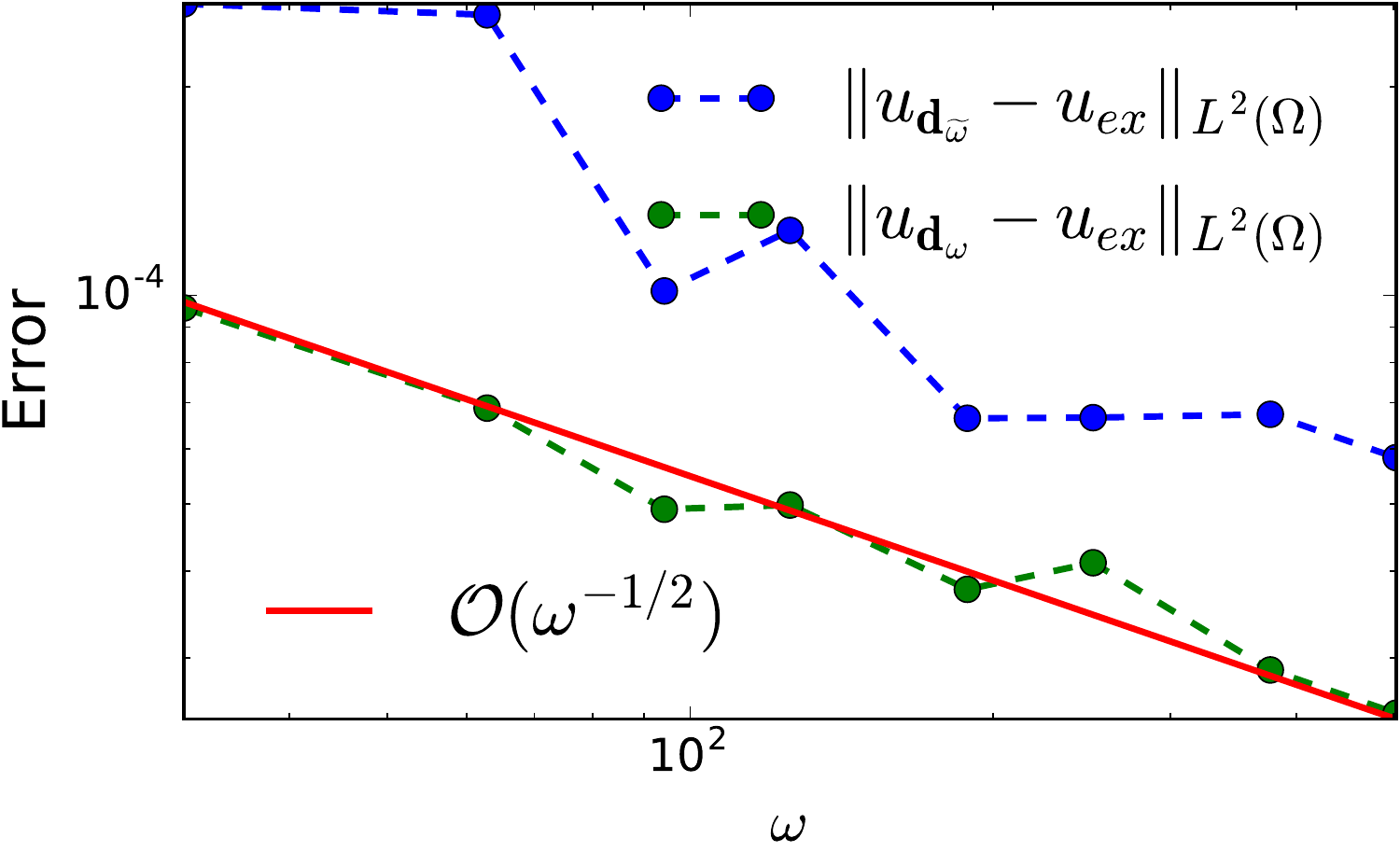}
\includegraphics[trim = 0mm 0mm 0mm 0mm,clip,width=7cm ]{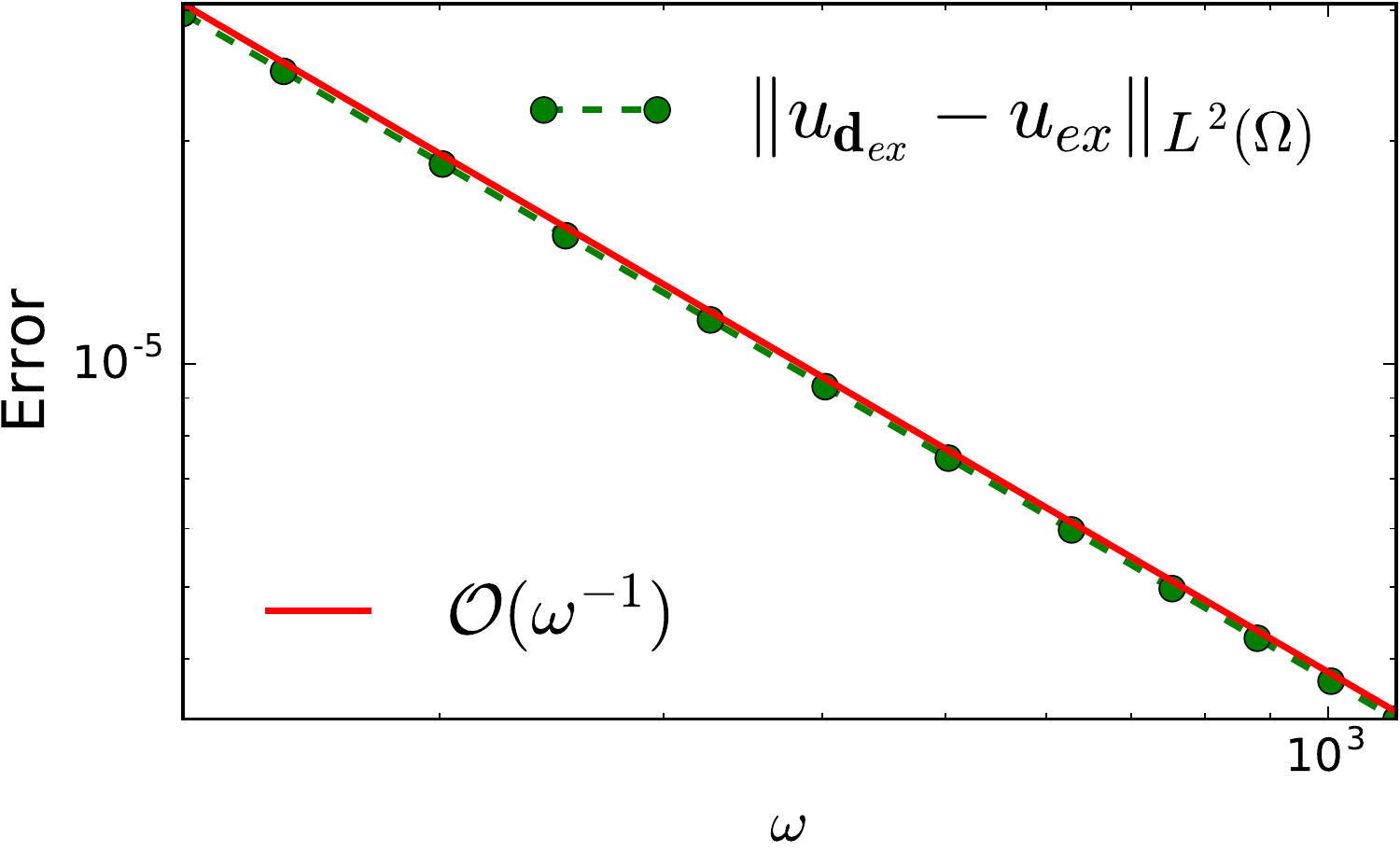}
\includegraphics[trim = 0mm 0mm 0mm 0mm,clip,width=7cm ]{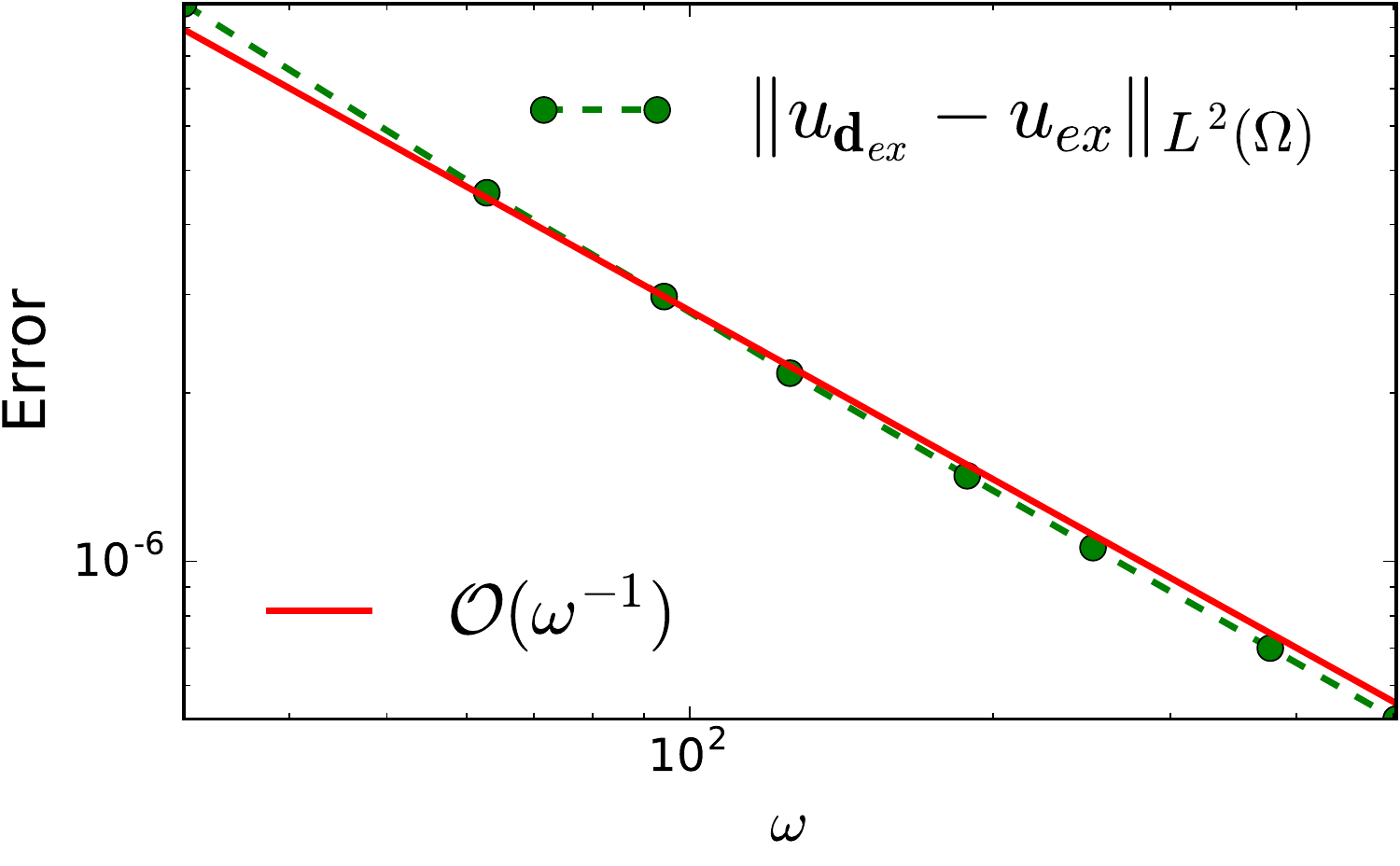}
}
\end{center}
\vspace{-0.5cm}
\caption{Tests with source outside domain, NPW = 6. Left: one point source; Right: four point sources. Top: ray direction errors; Middle: errors of ray-FEM solutions with ray directions estimated by NMLA; Bottom: errors of ray-FEM solutions with exact ray directions.}
\label{fig:source outside domain}
\end{figure}

% I'm not fully convinced about this...
Table \ref{tab:one source outside domain} and the left column of Figure \ref{fig:source outside domain} show that the asymptotic error, for both ray estimation and numerical solution by ray-FEM, decreases as the frequency increases. Moreover, they show that the ray-FEM algorithm is stable and quasi-optimal with fixed NPW, i.e., $\omega h = \cO (1)$.  % Here the Quasi-optimality is not exactly the same as Babuska's definition ??
Quasi-optimality holds if the ratio of numerical solution error to the best approximation error is bounded by a constant which is independent of frequency $\omega$. Here we assume the best approximation error in ray-FEM space $\inf_{u_h \in V^h_{Ray}(\mathcal{T}_h)} \| u - u_h \|$ is of the same order as the interpolation error $\| u - u_I^h \|$. So we can define the estimated quasi-optimality constant by $\frac{ \Vert u -  u_h \Vert_{L^2(\Omega)}}{\Vert u - u_I^h \Vert_{L^2(\Omega)} }$, and it is shown in Figure \ref{fig:ex1_opt_con}.
 Also it shows that one more iteration using iterative ray-FEM can significantly improve the final numerical solution to the order of  $\cO (\omega^{-1})$, which is of the same order when exact ray direction $\textbf{d}_{ex}$ is used in ray-FEM,  due to the asymptotic error for geometric ansatz. Note, however, that NMLA with curvature correction is only valid for a single point source in homogeneous media.

\begin{table}[h]
%title of the table
\centering % centering table
\renewcommand{\arraystretch}{1.2} % Default value: 1
\begin{tabular}{||c ||c|c|c|c|c||c} % creating eight columns
\hline\hline %inserting double-line
\
$\omega /2\pi$ & 20 & 40 & 80 & 160  \\
\hline
\
$1/h $  & 120 &  240  & 480  & 960 \\
\hline\hline      % inserts single-line  [0.5ex]
\
$\Vert \theta (\textbf{d}_{\widetilde{\omega}}) - \theta_{ex}\Vert_{L^2}$ &  7.50e-04  &  4.26e-04  &  1.96e-04  &  1.07e-04    \\
\hline
\
$\Vert\theta (\textbf{d}_{\omega}) - \theta_{ex}\Vert_{L^2}$ &  1.82e-04  &  7.99e-05  &  4.43e-05  &  2.10e-05   \\
[0.2ex] % [1ex] adds vertical space
\hline \hline
\
$\Vert u_{\textbf{d}_{\widetilde{\omega}}} - u_{ex}\Vert_{L^2}$  &  4.36e-05  &  1.92e-05  &  9.03e-06  &  4.69e-06 \\
\hline
\
$\Vert u_{\textbf{d}_{\omega}} - u_{ex} \Vert_{L^2}$&  3.15e-05  &  1.47e-05  &  7.57e-06  &  3.73e-06   \\
[0.2ex] % [1ex] adds vertical space
\hline \hline
\
$\Vert u_{\textbf{d}_{ex}} - u_{ex}\Vert_{L^2}$  &  2.97e-05  &  1.49e-05  &  7.47e-06  &  3.74e-06  \\
[0.2ex] % [1ex] adds vertical space
\hline \hline% inserts single-line
\end{tabular}
\caption{Errors of one point source problem for fixed NPW = 6. $\theta_{ex}$ is the exact ray angle, $\theta (\textbf{d}_{\widetilde{\omega}})$ and $\theta (\textbf{d}_{\omega})$ are ray angle estimations using low and high frequency waves, respectively;  $u_{\textbf{d}_{\widetilde{\omega}}}$, $u_{\textbf{d}_{\omega}} $ and $u_{\textbf{d}_{ex}} $ are ray-FEM solutions using low-frequency ray estimation $\textbf{d}_{\widetilde{\omega}}$, high frequency ray estimation $\textbf{d}_{\omega}$, and exact ray $\textbf{d}_{ex}$, respectively.}
\label{tab:one source outside domain}
\end{table}

\begin{figure}[h]
\centering
\includegraphics[scale=0.7]{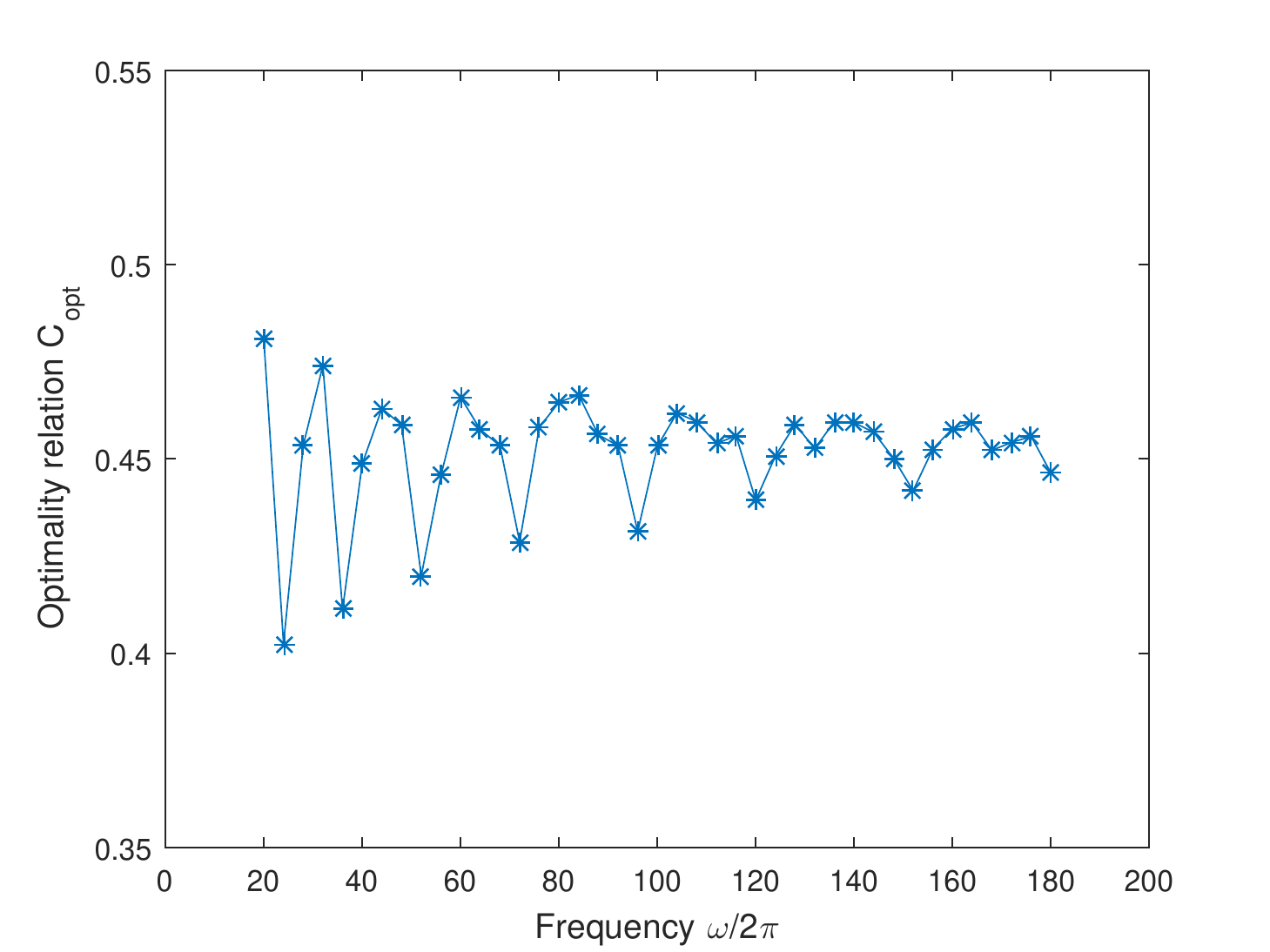}
\caption{The stars defined by $\frac{ \Vert u -  u_h \Vert_{L^2(\Omega)}}{\Vert u - u_I^h \Vert_{L^2(\Omega)} }$ with \mbox{NPW}= 6, give an indication of the optimality constant.}
\label{fig:ex1_opt_con}
\end{figure}

%(CASE 2 here???)

Next we show that our method can handle multiple wave fronts by probing the whole domain and extracting dominant ray directions locally.  The setup is exactly as above except that there are four point sources. The exact solution is given by
\begin{equation} \label{eq:exact_solution_4_point_source_outside}
\begin{array}{ll}
u_{ex}(x,y) &= \sqrt{\omega}H_0^{(1)}(\omega \sqrt{(x+20)^2 + (y+20)^2}) + 2\sqrt{\omega}H_0^{(1)}(\omega\sqrt{(x-20)^2 + (y-20)^2})\\
& + 0.5\sqrt{\omega}H_0^{(1)}(\omega\sqrt{(x+20)^2 + (y-20)^2}) - \sqrt{\omega}H_0^{(1)}(\omega\sqrt{(x-20)^2 + (y+20)^2}).
\end{array}
\end{equation}
The main difficulty of this example compared to the one above, is that the low-frequency wave solution by the standard FEM contains multiple wave fronts at each point due to the interference of multiple sources.
The numerical results are shown in right column of Figure \ref{fig:source outside domain}. In this case, NMLA with curvature correction does not apply so we use the the standard NMLA version for plane wave decomposition described in Section  \ref{subsection:NMLA} to estimate local dominant ray directions. As analyzed in Section \ref{subsec:numerical-ray-FEM} and \ref{appendix-b}, the expected error for ray direction estimation and numerical solution is of order $\cO (\omega^{-1/2})$ due to the curved wave fronts. The numerical results show that the ray-FEM meets the expectation without pollution as the frequency increases.

\subsection{Phase errors}
Here we show that by incorporating the estimated ray directions, ray-FEM can capture the phase much more accurately.
We test our algorithm with a source inside the domain. In particular, we use a point source, given its importance in many practical applications, in particular, in geophysics, in which the sources are often modeled as point sources. Moreover, in applications oriented towards inverse and imaging problems, having a numerical method that produces the correct phase in the far field is of great importance in order to properly locate features in the image.

% Computing the solution to the Helmholtz equation using point sources is of great importance in many field (which fields? refereces? ), given that by simplicity many different sources are modeled using point sources. Moreover, for practical applications, it is extremely important to obtain the phases in the far field correctly.

\begin{figure}[h]
\begin{center}
{\includegraphics[trim = 0mm 0mm 0mm 0mm,clip,width=7cm ]{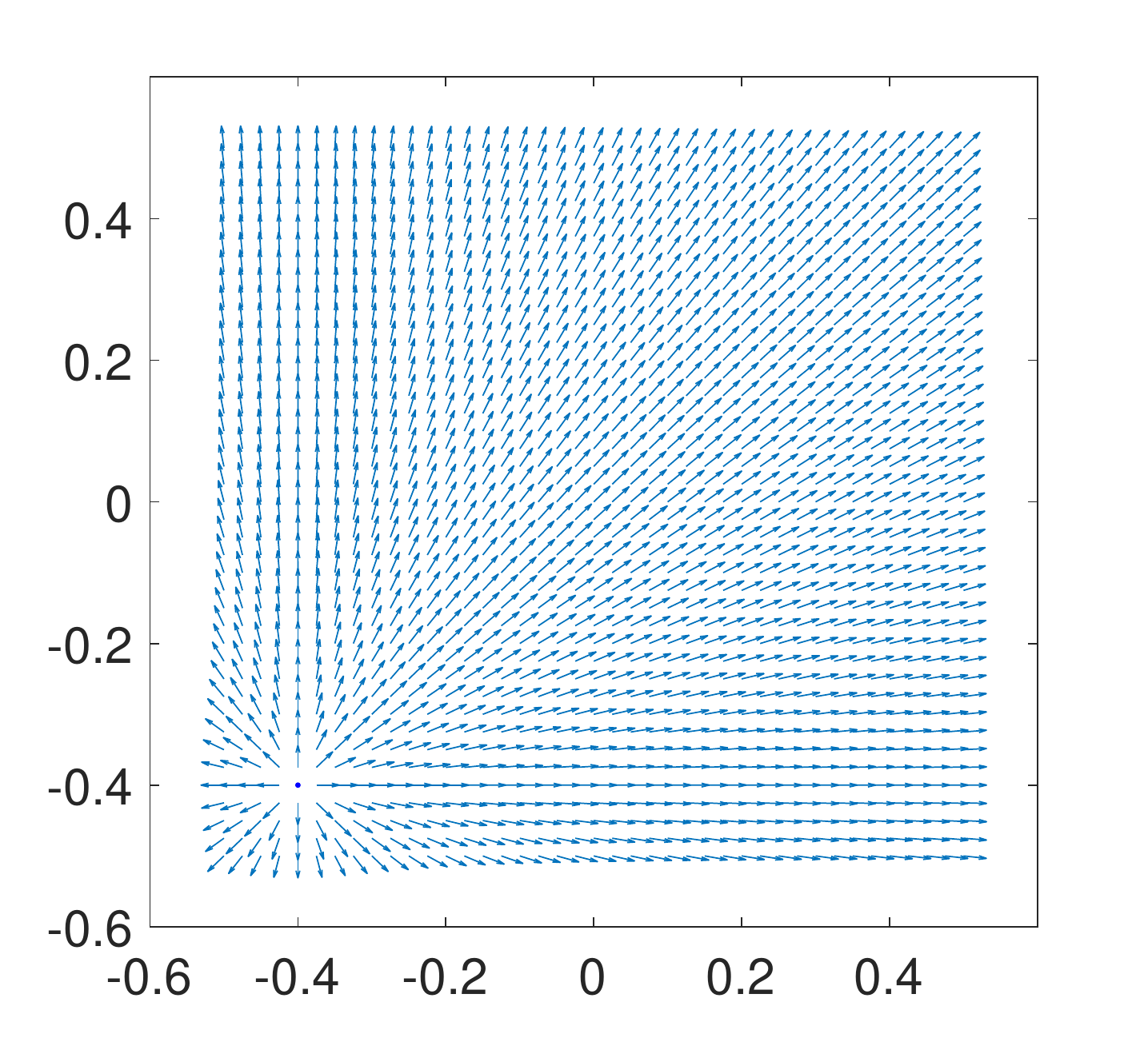}
\includegraphics[trim = 0mm 0mm 0mm 0mm,clip,width=7.45cm ]{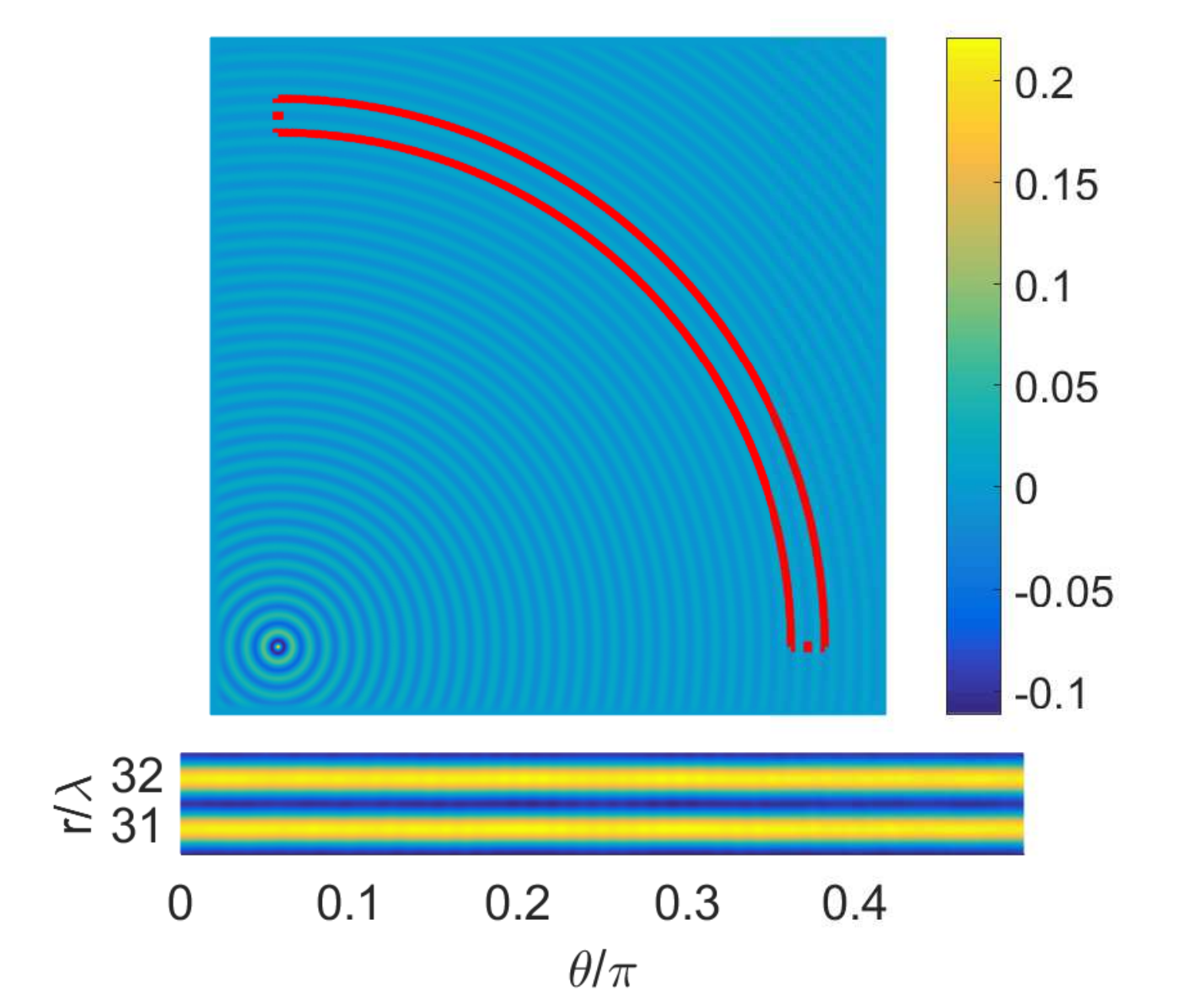}
 }
\end{center}
\vspace{-0.5cm}
\caption{One point source inside homogeneous medium domain, $\omega = 80\pi$, NPW = 6. Left: ray direction field captured by NMLA; Right: polar plot of ray-FEM solution, $r/ \lambda$: the number of wavelength away from the source.} \label{fig:one point source inside homogeneous domain: 80pi}
\end{figure}

In this experiments we focus our attention on the far-field since our current method can not deal with singularities in amplitude and phase at source points. We start by solving the Helmholtz equation using a slight modification of Algorithm \ref{alg:Iter_Ray_FEM} and using standard finite elements at the source point. In this case the source term is located inside the domain $\textbf{f }= \delta (\textbf{x} - (-0.4,-0.4))$, but we will use the associated column of the mass matrix (normalized by mesh size $h$) as the right hand side. For vertices near the source, we use the exact ray direction, the radial direction, in our ray-FEM.
And we find the ray directions by NMLA for vertices away from the source; see Figure \ref{fig:one point source inside homogeneous domain: 80pi} left part for the ray direction field.

To demonstrate the phase errors in numerical solutions, we plot the computed wave field on a 90 degree part of an annulus \cite{Stolk}, with the radial coordinate varying on an interval of about two wavelengths; see Figure \ref{fig:one point source inside homogeneous domain: 80pi} right part.

\begin{figure}[h]
\centering
\includegraphics[scale=1.0]{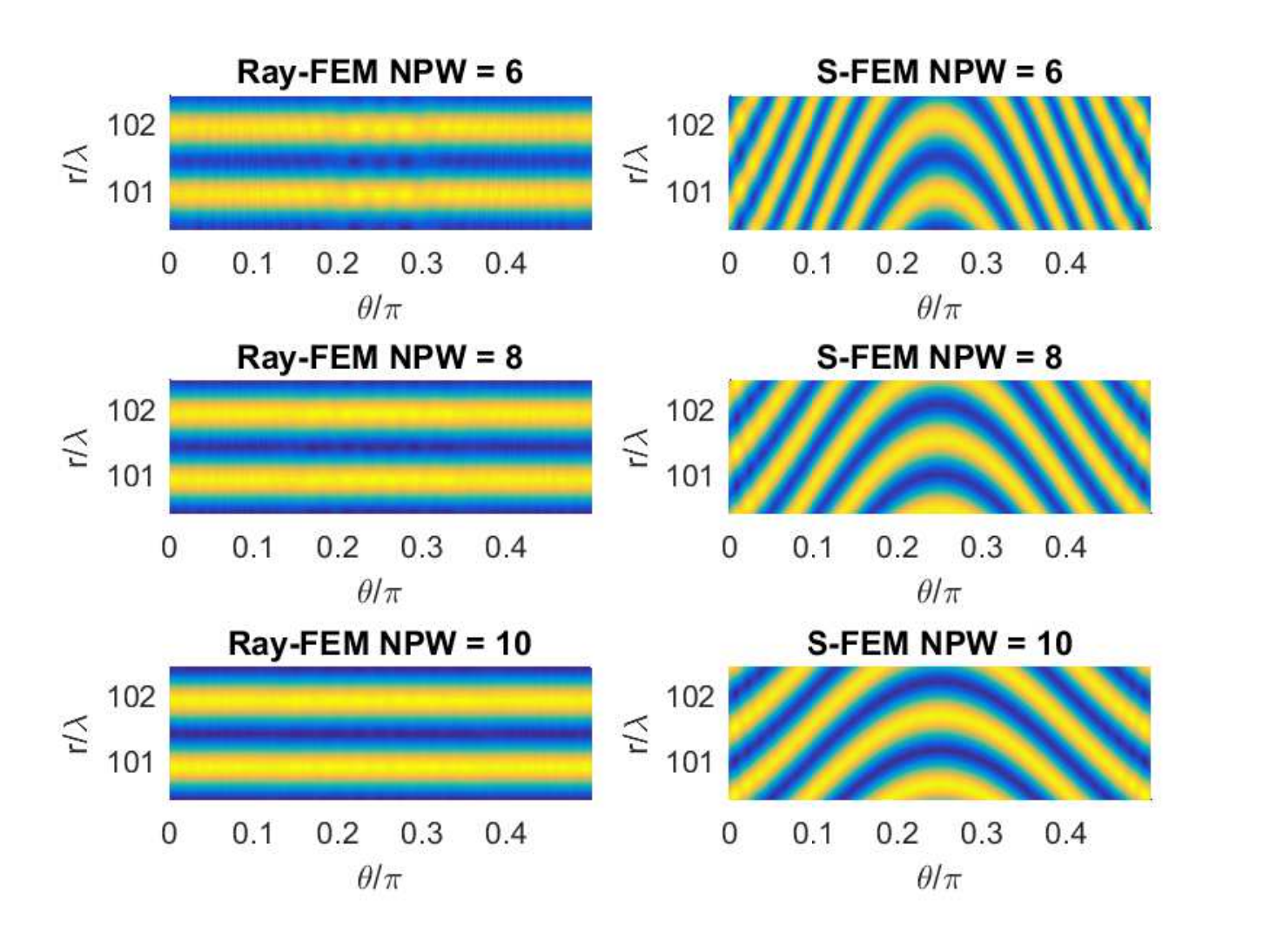}
\vspace{-0.5cm}
\caption{Polar plot of ray-FEM solution $\omega = 250\pi$. $r/ \lambda$: the number of wavelength away from the source.  }
\label{fig:one point source inside homogeneous domain: 250pi}\end{figure}

In this case the frequency $\omega = 250 \pi$ is fixed, but we increase the number of grid points per wavelength. Figure \ref{fig:one point source inside homogeneous domain: 250pi} depicts the behavior of both ray-FEM solution and standard FEM solution. From the figure we can easily observe the superiority of the ray-FEM on minimizing the phase error, even using relatively coarse meshes.

%\vspace{20mm}

% the resulting wavefiled on a 90 degree part of an annulus (red box in Fig (show Fig)), with radial coordinate varying on an interval of about two wavelengths. The point source locates at $(-0.3,-0.3)$. Fig blah depicts the behavior of both Ray based finite elements and standard one. From fig we can easily observe the superiority of the ray based FEM on minimizing the phase error, even using relatively coarse meshes.

\begin{figure}[h]
\begin{center}
{\includegraphics[trim = 0mm 0mm 0mm 0mm,clip,width=7cm ]{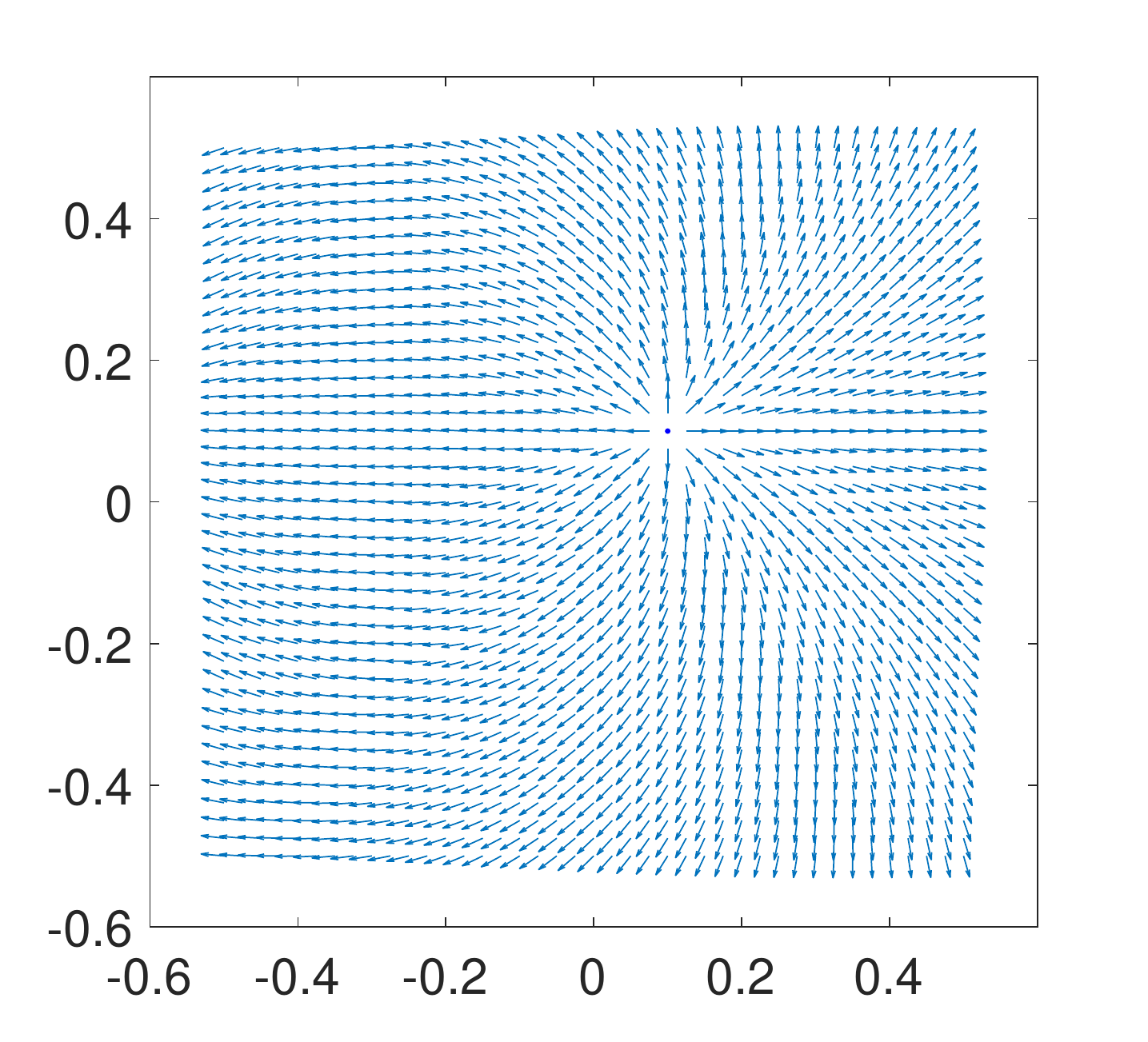}
\includegraphics[trim = 0mm 0mm 0mm 0mm,clip,width=7.5cm ]{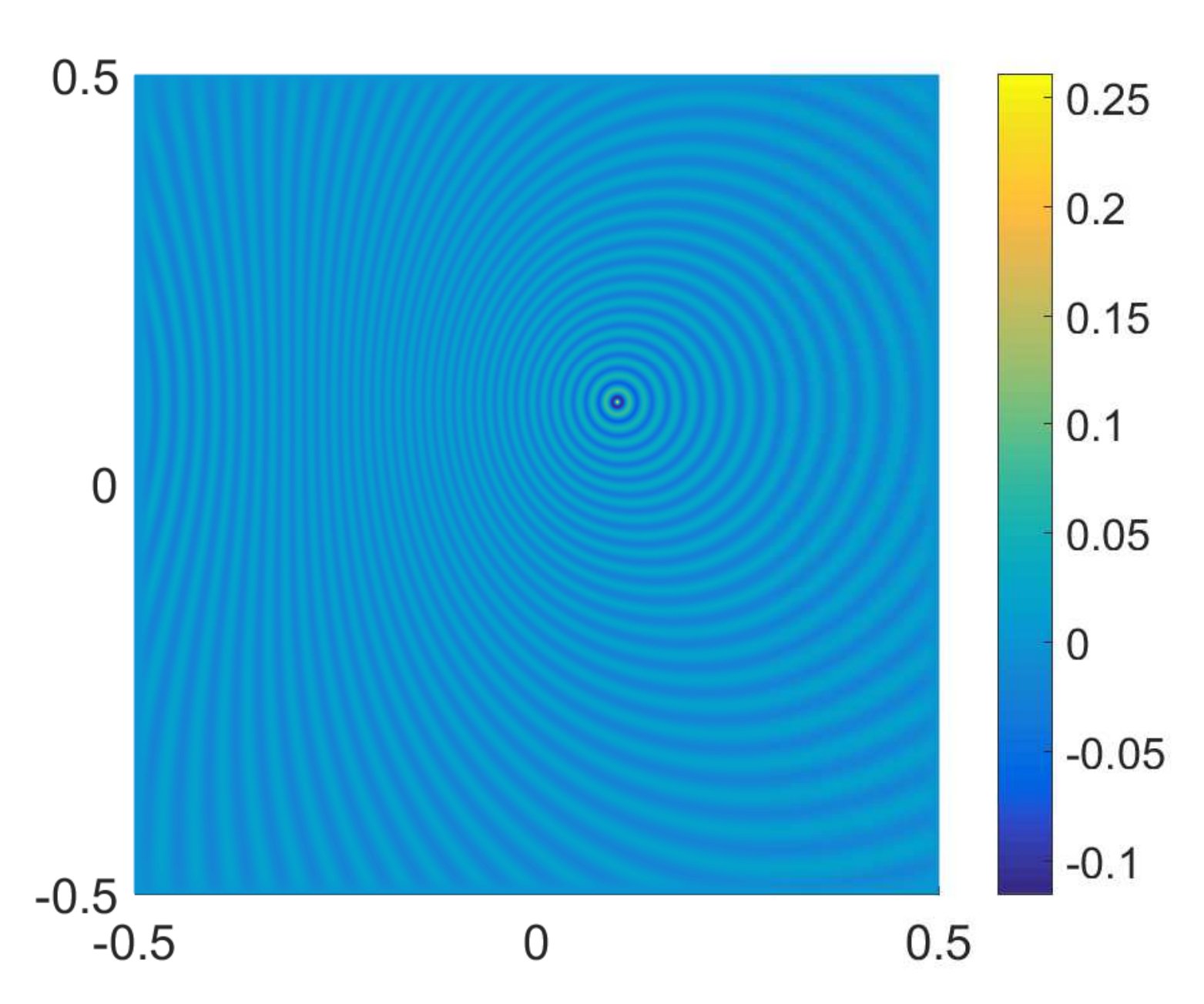}
}
\end{center}
\vspace{-0.5cm}
\caption{One point source inside heterogeneous medium domain with Gaussian wave speed $c(x,y) = 3 - 2.5e^{ -((x +0.125)^2 + (y-0.1)^2)/0.8^2}$,  $\omega = 80\pi$, NPW = 10. Left: ray direction field captured by NMLA;  Right: wave field computed by ray-FEM.}
\label{fig:one point source inside heterogeneous domain: 80pi}
\end{figure}

\begin{figure}[h]
\begin{center}
{\includegraphics[trim = 0mm 0mm 0mm 0mm,clip,width=7.2cm ]{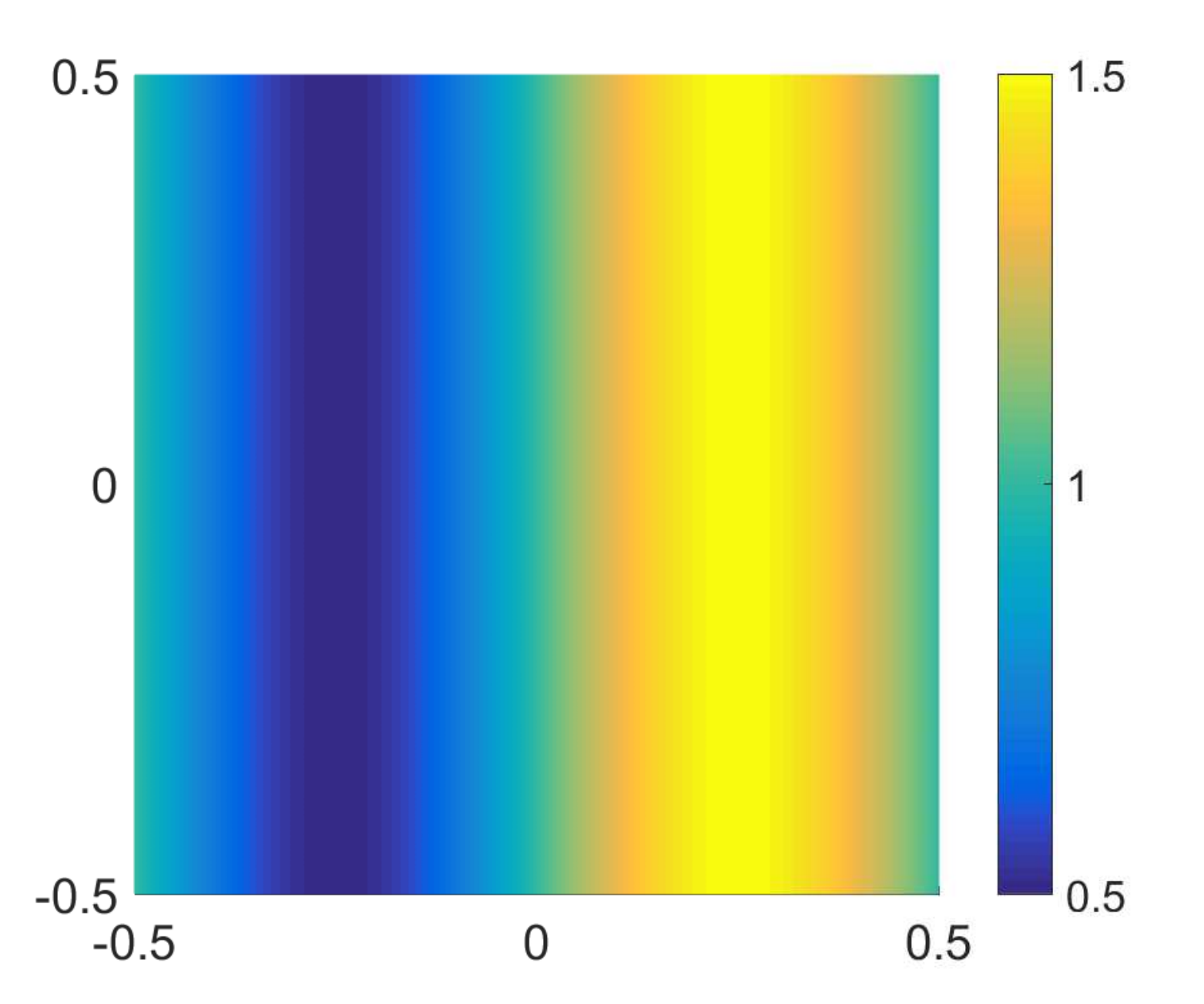}
\includegraphics[trim = 0mm 0mm 0mm 0mm,clip,width=7.2cm ]{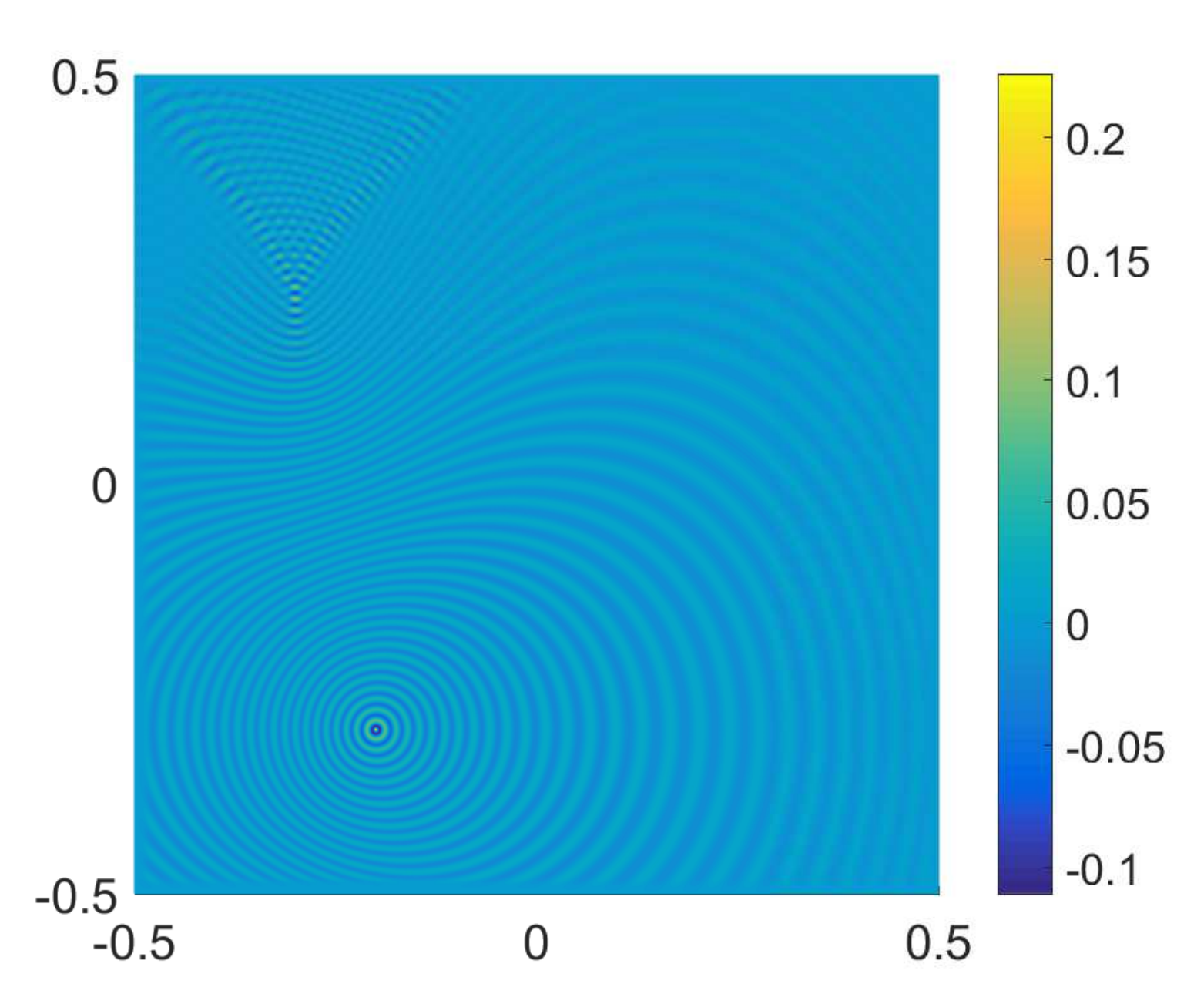}
}
\end{center}
\vspace{-0.5cm}
\caption{One point source inside heterogeneous medium domain with sinusoidal wave speed $c(x,y) = 1+ 0.5\sin(2\pi x)$,  $\omega = 80\pi$, NPW = 10. Left: wave speed; Right: wave field computed by ray-FEM.}
\label{fig:caustics: 80pi}
\end{figure}

We then solve the Helmholtz equation using a heterogeneous medium given by Figure \ref{fig:one point source inside heterogeneous domain: 80pi} with the source located inside. We also provide an experiment where we show the ability of the method in this paper to handle wave field with caustics; see Figure \ref{fig:caustics: 80pi}. Again radial directions are used for local ray directions near the source point.

\iffalse
\begin{figure}[h]
\centering
\includegraphics[scale=0.7]{homogeneous_ray_field.pdf}
%\includegraphics[scale=0.7]{eps/homogeneous_wave_field.pdf}
%\includegraphics[scale=0.7]{eps/caustics_wave_field.pdf}
%\includegraphics[scale=0.7]{eps/caustics_wave_speed.pdf}
\caption{The stars give the optimality constant $C_{opt} = \frac{ \Vert u -  u_h \Vert_{L^2(\Omega)}}{\Vert u - u_I \Vert_{L^2(\Omega)} }$ with \mbox{NPW}= 6. }
\label{fig:ex1_opt_con}
\end{figure}
\fi

\subsection{Complexity tests} \label{section:experiments_fast_methods}
In this subsection we test the computational complexity for ray-FEM.  A key step is solving the sparse linear systems generated by ray-FEM using iterative methods with a performant preconditioner, e.g., domain decomposition techniques coupled with high-quality absorbing/transmission boundary conditions. In our tests, we use a modification of the method of polarized traces to solve the linear systems resulting from both standard FEM and ray-FEM as described in Section \ref{subsection:fast_linear_solver}.

We use the numerical experiments to demonstrate the overall computational complexity of our ray-FEM method.
In particular, we solve the Helmholtz equation with a point source in both homogeneous medium and heterogeneous medium. We compute for many different frequencies, using Algorithm \ref{alg:Iter_Ray_FEM} with only one iteration of ray-FEM, the solution to the Helmholtz equation posed on $\Omega = [-0.5, 0.5] \times [-0.5, 0.5]$ with absorbing boundary conditions implemented via PML. For each frequency we report the execution time of the low and high frequency problems and the time spent in processing the data using NMLA to extract the dominant ray information.

As explained in Section \ref{Section:Algorithms}, in order to process the data using NMLA we need to solve the low-frequency problem in a slightly larger domain. The size of the larger domain is given by the sampling radius of the NMLA. For the sake of simplicity, we use a low-frequency subdomain, $\Omega_{low} = [-1,1]\times [-1,1]$, i.e., four times bigger than the original domain. The size can be reduced in order to lower computational cost for the low-frequency problem.

\begin{figure}[h]
\begin{center}
{\includegraphics[trim = 0mm 0mm 0mm 0mm,clip,width=7cm ]{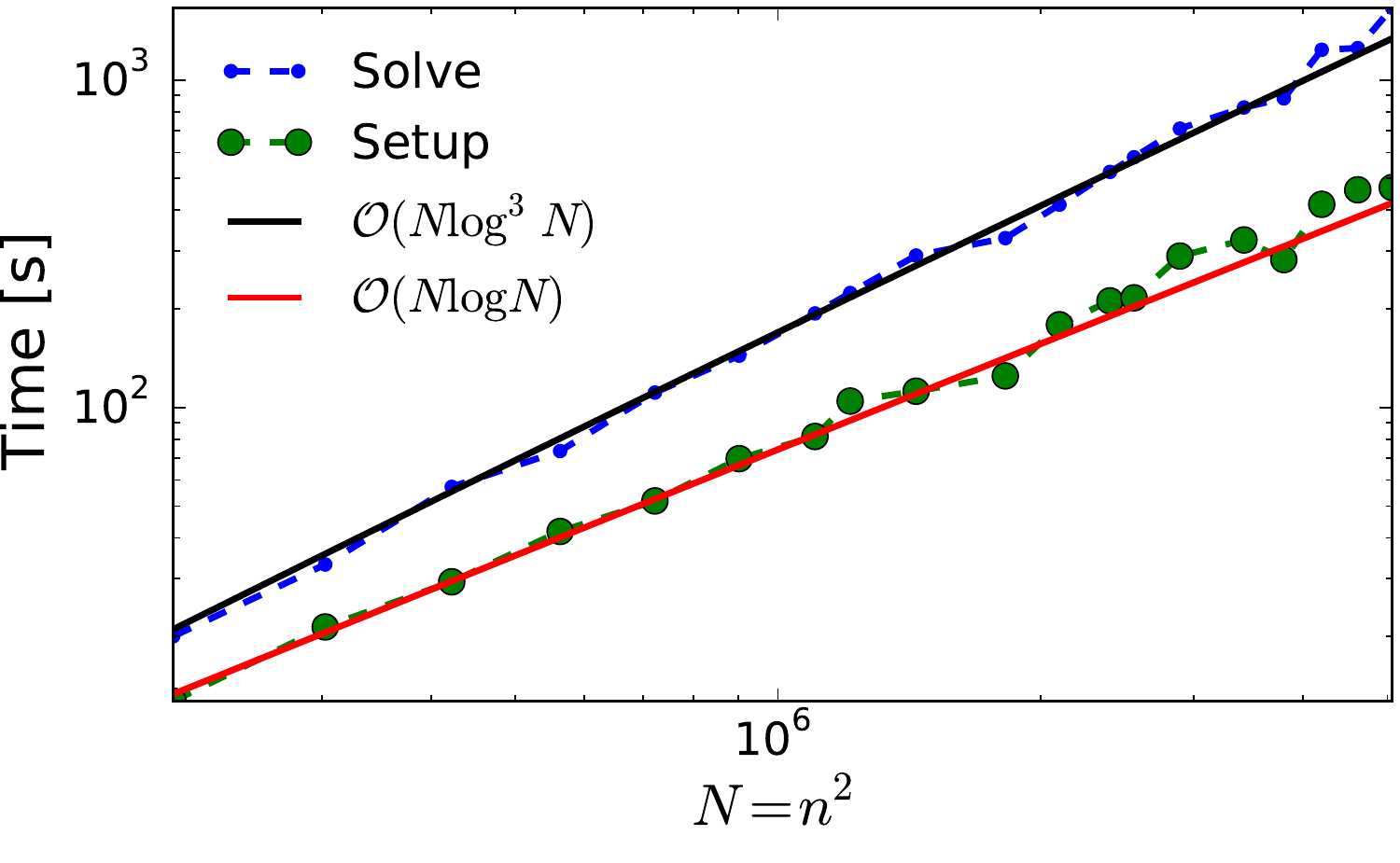}
\includegraphics[trim = 0mm 0mm 0mm 0mm,clip,width=7cm ]{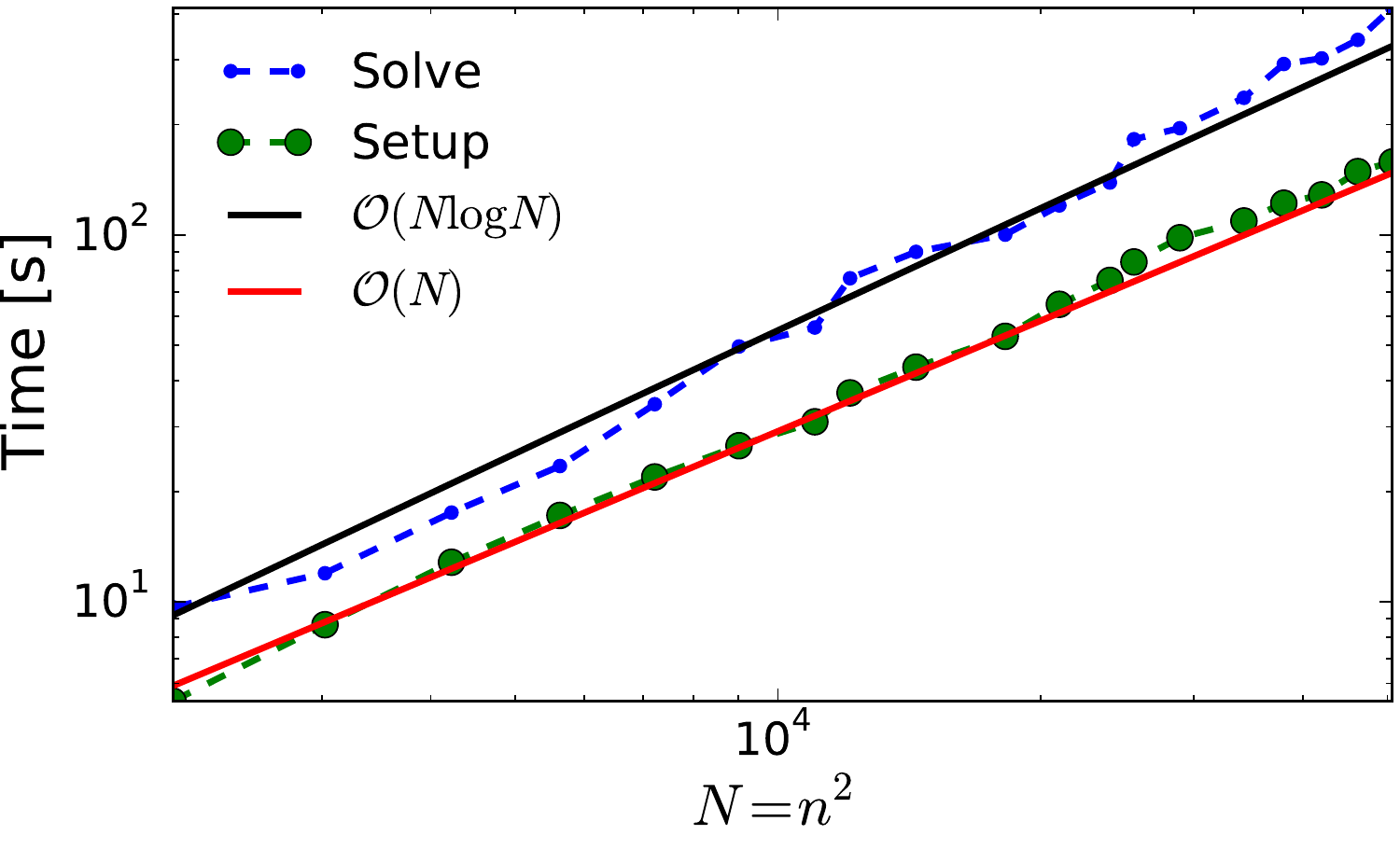} }
\end{center}
\vspace{-0.5cm}
\caption{Runtime for solving the Helmholtz equation with homogeneous wave-speed using GMRES preconditioned with the the method of polarized traces. The tolerance was set up to $10^{-7}$. Left: runtime for solving the low-frequency problem. Right: Runtime for solving the high-frequency problem with the adaptive basis.} \label{fig:fast_solver_homogenenous}
\end{figure}

\begin{figure}[h]
\begin{center}
{\includegraphics[trim = 0mm 0mm 0mm 0mm,clip,width=7cm ]{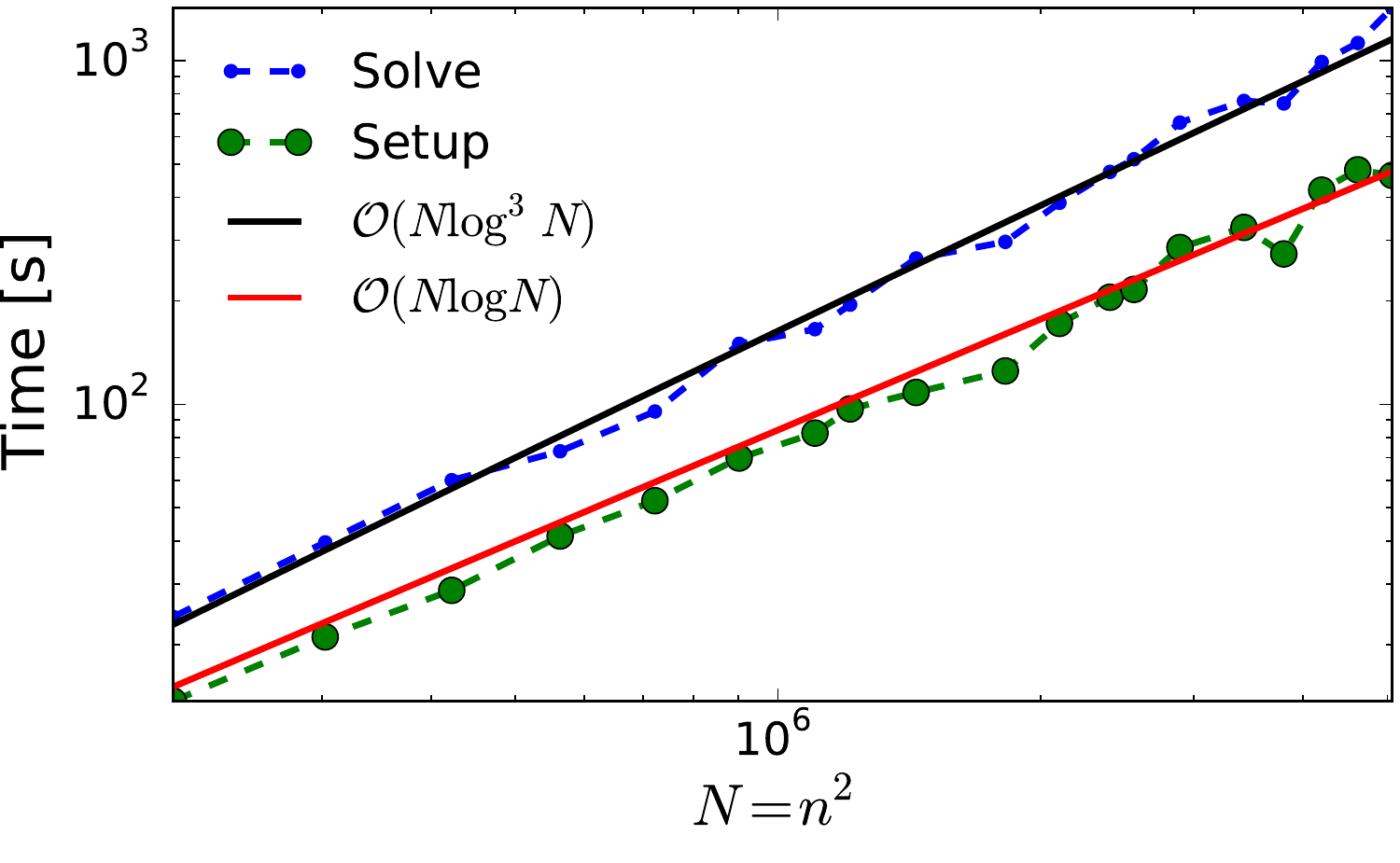}
 \includegraphics[trim = 0mm 0mm 0mm 0mm,clip,width=7cm ]{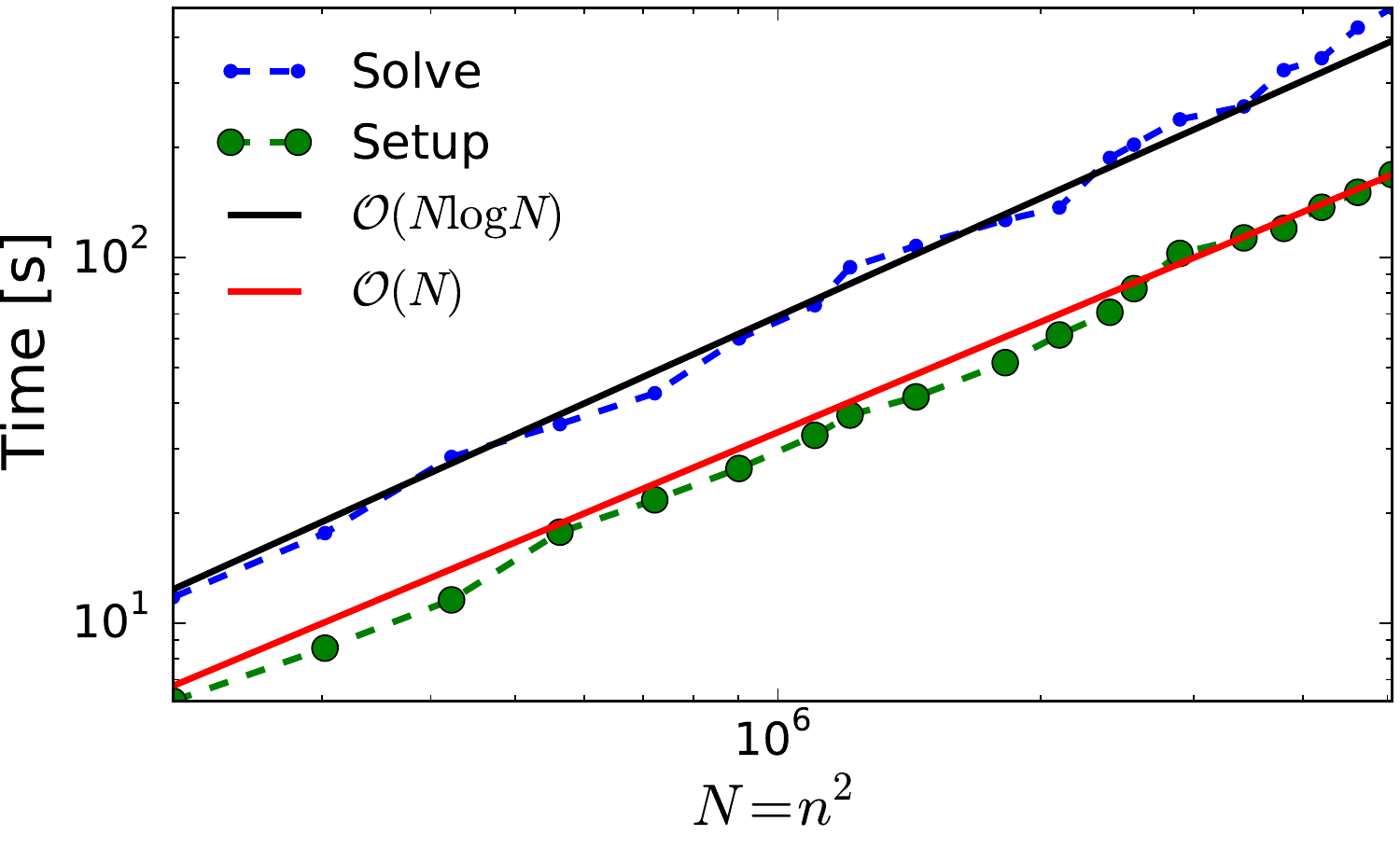}
}
\end{center}
\vspace{-0.5cm}
\caption{Runtime for solving the Helmholtz equation with heterogeneous wave-speed using GMRES preconditioned with the the method of polarized traces. The tolerance was set up to $10^{-7}$. Left: runtime for solving the low-frequency problem. Right: runtime for solving the high-frequency problem with the adaptive basis.} \label{fig:fast_solver_heterogeneous}
\end{figure}

The main issue with the low-frequency solver in our case are the PML's, given that each thin slab contains less than a wavelength across, the PML may not be very effective. In order to decrease the number of iterations to converge, we increase the PML points logarithmically with the frequency. This implies a slightly more expensive setup cost and solve cost as shown in Figures~\ref{fig:fast_solver_homogenenous} left and \ref{fig:fast_solver_heterogeneous} left.

Figure~\ref{fig:fast_solver_homogenenous} shows the runtime for solving the Helmholtz equation with a point source inside a homogeneous medium. We can observe that the overall cost is $\cO(N)$ up to poly-logarithmic factors as shown in our complexity study. The low-frequency solver has a slightly higher asymptotic cost in this case, given the ratio between the width of the PML and the characteristic wavelength inside the domain.

Figure~\ref{fig:fast_solver_heterogeneous} shows the runtime for solving the Helmholtz equation with a point source inside a heterogenous medium. We can observe the same scaling as before, albeit with slightly larger constants.

\clearpage

\section{Conclusion} \label{Section:Conclusion}

In this work we present a numerical method, the ray-FEM, for the high frequency Helmholtz equation in smooth media based on learning problem specific basis functions to represent the wave field. The key information, local ray directions, is extracted from a relative low frequency wave field that has probed the whole domain. These local ray directions are then incorporated into the basis to improve both stability and accuracy in the computation for high frequency wave field.  Moreover, both local ray directions and the high frequency wave field can be further improved through more iterations. Numerical tests suggest that our method only requires a fixed number of points per wave length without pollution effect as frequency becomes large. By designing a fast solver for the discretized linear systems an overall complexity of order $\cO(\omega^d\log \omega)$ is achieved.

However, our ray-FEM can not handle singularities of both the amplitude and phase on a mesh. We will develop a hybrid method that combines local asymptotic expansion near the source and the ray-FEM away from the source in our future work.

\clearpage

%\nocite{*}

\section*{Acknowledgments}
Zhao is partially supported by NSF grant (1418422). Qian is partially supported by NSF grants (1522249 and 1614566).




%


\section{References}
\bibliographystyle{plain}
\bibliography{references}

%!TEX root = draft.tex

\appendix

%\section*{Appendix}

\section{Stability and error analysis for NMLA}\label{appendix-a}
In this section we summarize the stability result and error estimate from \cite{Benamou:NMLA_revisited} for completeness.
For simplicity we use the single wave case, i.e., $N = 1$. Moreover, we assume the measurement data is a perturbation to the perfect plane wave data of the form $U(\theta) =  U^{plane}(\theta) + \delta U(\theta)$, where $U^{plane}$ denotes a single plane wave data in the form of \eqref{eq:model_NMLA}. Let $\theta^*$ denote the angle for which $\theta \mapsto \mathcal{B}U(\theta)$ is maximum. Assuming that the noise level satisfies
\begin{equation} \label{eq:NMLA_error_condition}
||\delta U||_{L^{\infty}} < \frac{1}{4B^*}|B_1|,
\end{equation}
where $B^*\leq 0.89$ is a pure constant independent of $\omega$ and $B_1$ is the complex amplitude of the plane wave. Then the error in the angle estimation is given by
\begin{equation}
|\theta_1 - \theta^*| \leq \frac{2\pi}{2L_{\alpha}+1} \sim \mathcal{O}(\frac{1}{\alpha}),
\quad \alpha = kr \sim \infty.
\end{equation}

Similar results can be derived for multiple waves $N > 1$. We remark that $\frac{1}{4B^*} \simeq 0.28$, which implies that if the relative noise level does not surpass $28 \%$ the angle will be detected within an error of order $\cO(\frac{1}{kr})$. In Benamou's work \cite{Benamou:NMLA_revisited}, an analysis of a point source shows that $|\theta_1 - \theta^*|$ decreases like $\mathcal{O}(\omega^{-1/2})$ when the point $\x_0$ is far away from the source and the radius of the observation circle is chosen like $r \sim \omega^{-1/2}$ for large $\omega$. We obtain similar accuracy order for general noisy plane waves under some smoothness conditions, see details in \ref{appendix-b}.

\section{Error analysis of wave-field as a perturbed plane wave data} \label{appendix-b}
As introduced in Section \ref{subsection:NMLA}, NMLA is a tool to process a signal that is (approximately) a superposition of plane waves with frequency $\omega$ and to extract each plan wave component by sampling the signal on a circle/sphere with radius $r$ around a reference point. As shown in \ref{appendix-a}, provided that the perturbation of the signal is relative small compared to the signal, the estimation of the plane wave directions converges and the error is $\cO(\frac{1}{\omega r})$. In this application, we use NMLA to process wave-field data, which is the numerical solution to the Helmholtz equation, to extract the directions of dominant wave fronts based on the geometric optics ansatz \eqref{eq:plane_wave_expansion} in the high-frequency regime. Hence it is important to study the wave field data as a perturbation of plane wave data locally and estimate the error in the ray directions obtained from NMLA. In particular, this analysis allows us to find the optimal choice of the radius of the sampling circle/sphere, in order to achieve the minimal asymptotic error for the ray direction estimation in terms of the frequency $\omega$ of the Helmholtz equation which generates the wave-field data. The result is crucial for both error analysis and implementation of ray-FEM.  Since the wave field data in our application is the numerical solution to the Helmholtz equation, its perturbation can be composed as the sum of three components:
\begin{enumerate}
\item numerical error in solving the Helmholtz equation and interpolation error in obtaining data on the sampling circle/sphere for NMLA from the numerical solution on a fixed mesh,
\item the asymptotic error in the geometric optics ansatz,
\item the local deviation of a smooth curved wave front from a planar wave front.
\end{enumerate}

On a mesh with mesh size $h=\cO(\omega^{-1})$, the last components, which we call the phase error, is the dominant factor among the three. We present below an analysis of the phase error, in which, for simplicity, we only consider one wave front.

Let consider a single wave front, $u(\textbf{x}) = A(\textbf{x})e^{i\omega \phi(\textbf{x})}$; following the notation used throughout the paper, assume the reference point to be $\x_0$, and the small sampling circle around $\x_0$ to be $\{\x|\textbf{x} - \textbf{x}_0 = r \widehat{\textbf{s}}\}$, $\nabla \phi (\textbf{x}_0) = \eta_0 \widehat{\textbf{d}}_0$, where $r\ll 1, |\widehat{\textbf{s}}| = 1$,  $\eta_0 = 1/c(\textbf{x}_0)$, $|\widehat{\textbf{d}}_0 |= 1$.
\begin{displaymath}
\vspace{2mm}
\begin{array}{ll}
A(\textbf{x}) & =
A(\textbf{x}_{0}) +  \nabla A( \textbf{x}_0) \cdot \left( \textbf{x} - \textbf{x}_0 \right) + O \left( (\textbf{x}-\textbf{x}_0)^2 \right)
=
A(\textbf{x}_{0}) + r \left( \nabla A( \textbf{x}_0) \cdot  \widehat{\textbf{s}} \right) + O \left( r^2 \right),
\vspace{2mm}
\\
\vspace{2mm}
\phi(\textbf{x})  &= \phi(\textbf{x}_{0}) + \nabla \phi (\textbf{x}_0) \cdot (\textbf{x} - \textbf{x}_0) + \frac{1}{2} \left(\textbf{x} - \textbf{x}_0 \right)^T \nabla^2 \phi ( \textbf{x}_0)\left(\textbf{x} - \textbf{x}_0  \right) + O \left( (\textbf{x}-\textbf{x}_0)^3 \right)\\
\vspace{2mm}
&= \phi(\textbf{x}_{0}) + r \eta_0 \left(  \widehat{\textbf{d}}_0 \cdot  \widehat{\textbf{s}} \right) + \frac{1}{2}r^2 \left(\widehat{\textbf{s}}^T \nabla^2 \phi ( \textbf{x}_0) \widehat{\textbf{s}}  \right)  +O \left( r^3 \right).
\end{array}
\end{displaymath}
Denote $\phi_{0}(\textbf{x}) =  \phi(\textbf{x}_{0}) + \nabla \phi (\textbf{x}_0) \cdot (\textbf{x} - \textbf{x}_0) $, $u_0 (\textbf{x}) = A(\textbf{x}_{0}) e^{i \omega \phi_{0}(\textbf{x})}$, we have
\begin{displaymath}
\vspace{2mm}
\begin{array}{ll}
 \delta u(\textbf{x}) & = u(\textbf{x}) - u_0(\textbf{x})
\vspace{2mm}
\\
& = A(\textbf{x})e^{i\omega \phi(\textbf{x})}  - A(\textbf{x}_0)e^{i\omega  \phi_0 (\textbf{x})}
\vspace{2mm}
\\
& = \left[ A(\textbf{x}_0)e^{i\omega \phi(\textbf{x})} + r \left( \nabla A( \textbf{x}_0) \cdot  \widehat{\textbf{s}} \right) e^{i\omega  \phi(\textbf{x})} + O \left( r^2 \right)  \right] - A(\textbf{x}_0)e^{i\omega  \phi_{0}(\textbf{x})}
\vspace{2mm}
\\
& = A(\textbf{x}_0)e^{i\omega \phi_{0}(\textbf{x})} \left(  e^{ i\omega \left[  \frac{1}{2}r^2 \left( \widehat{\textbf{s}}^T \nabla^2 \phi ( \textbf{x}_0) \widehat{\textbf{s}}  \right)  +O \left( r^3 \right) \right] } - 1 \right)
 + r \left( \nabla A( \textbf{x}_0) \cdot  \widehat{\textbf{s}} \right) e^{i\omega \phi(\textbf{x})} + O \left( r^2 \right),
 \end{array}
 \vspace{2mm}
\end{displaymath}
\begin{displaymath}
\begin{array}{ll}
\frac{\partial}{\partial r} \left( \delta u(\textbf{x})\right) & =
\frac{\partial}{\partial r} \left( A(\textbf{x})e^{i\omega \phi(\textbf{x})}  - A(\textbf{x}_0)e^{i\omega  \phi_0 (\textbf{x})} \right)
\vspace{2mm}
\\
& = \left( \nabla A(\textbf{x}_0) \cdot \widehat{\textbf{s}} + O(r) \right)e^{i\omega \phi (\textbf{x})} + A(\textbf{x})e^{i\omega \phi (\textbf{x})}i\omega \left[ \eta _0 (\widehat{\textbf{d}}_0 \cdot \widehat{\textbf{s}}) + r \left( \widehat{\textbf{s}}^T \nabla ^2 \phi (\textbf{x}_0) \widehat{\textbf{s}}  \right)  + O(r^2) \right]
\vspace{2mm}
\\
& \qquad  - A(\textbf{x}_0)e^{i\omega \phi_{0}(\textbf{x}_0) } i\omega \eta _0 (\widehat{\textbf{d}}_0 \cdot \widehat{\textbf{s}})

\vspace{2mm}
\\
&= \left( \nabla A(\textbf{x}_0) \cdot \widehat{\textbf{s}} + O(r) \right)e^{i\omega \phi (\textbf{x})} + A(\textbf{x})e^{i\omega \phi (\textbf{x})}i\omega \left[ r \left(  \widehat{\textbf{s}}^T \nabla ^2 \phi (\textbf{x}_0)  \widehat{\textbf{s}}  \right)  + O(r^2) \right]
\vspace{2mm}
\\
& \qquad + \left( A(\textbf{x})e^{i\omega \phi(\textbf{x})} - A(\textbf{x}_0)e^{i\omega \phi_{0}(\textbf{x}_0)} \right)i\omega \eta _0 ( \widehat{\textbf{d}}_0 \cdot  \widehat{\textbf{s}})

\vspace{2mm}
\\
&= \left( \nabla A(\textbf{x}_0) \cdot  \widehat{\textbf{s}} + O(r) \right)e^{i\omega \phi (\textbf{x})} + A(\textbf{x})e^{i\omega \phi (\textbf{x})}i\omega \left[ r \left(  \widehat{\textbf{s}}^T \nabla ^2 \phi (\textbf{x}_0)  \widehat{\textbf{s}}  \right)  + O(r^2) \right]
\vspace{2mm}
\\
& \qquad + i\omega \eta _0 ( \widehat{\textbf{d}}_0 \cdot  \widehat{\textbf{s}}) \delta u(\textbf{x}).

\end{array}
\end{displaymath}
Then
\begin{displaymath}
\begin{array}{ll}
\vspace{2mm}
\delta U (\textbf{x})
& = \left( \frac{1}{i\omega \eta_0}\frac{\partial}{\partial r}  + 1 \right)  \delta u(\textbf{x})  \\

\vspace{3mm}
& = \frac{1}{i\omega \eta_0} \left( \nabla A(\textbf{x}_0) \cdot  \widehat{\textbf{s}} + O(r) \right)e^{i\omega  \phi (\textbf{x})}
+ \frac{1}{\eta_0} A(\textbf{x})e^{i\omega  \phi (\textbf{x})} \left[ r \left(  \widehat{\textbf{s}}^T \nabla ^2 \phi (\textbf{x}_0)  \widehat{\textbf{s}}  \right)  + O(r^2) \right]
\\
\vspace{3mm}
& \qquad + ( \widehat{\textbf{d}}_0 \cdot \textbf{s} ) \delta u(\textbf{x}) + \delta u(\textbf{x})

\\
\vspace{3mm}
& =  \frac{1}{i\omega \eta_0} \left( \nabla A(\textbf{x}_0) \cdot  \widehat{\textbf{s}} + O(r) \right)e^{i\omega  \phi (\textbf{x})}
+ \frac{1}{\eta_0} A(\textbf{x})e^{i\omega  \phi (\textbf{x})} \left[ r \left(  \widehat{\textbf{s}}^T \nabla ^2 \phi (\textbf{x}_0)  \widehat{\textbf{s}}  \right)  + O(r^2) \right]
\\
\vspace{3mm}
& \qquad + ( \widehat{\textbf{d}}_0 \cdot \widehat{\textbf{s}} + 1)  \left \{  A(\textbf{x}_0)e^{i\omega  \phi_{0}(\textbf{x})} \left(  e^{ i\omega  \left[  \frac{1}{2}r^2 \left(  \widehat{\textbf{s}}^T \nabla^2 \phi ( \textbf{x}_0)  \widehat{\textbf{s}}  \right)  +O \left( r^3 \right) \right] } - 1 \right)  \right.
\\
\vspace{3mm}
& \qquad \left. + r \left( \nabla A( \textbf{x}_0) \cdot   \widehat{\textbf{s}} \right) e^{i\omega \phi(\textbf{x})} + O \left( r^2 \right) \right \}.
\end{array}
\end{displaymath}
Hence
\begin{equation}
\vspace{3mm}
\begin{array}{ll}
|\delta U(\textbf{x})|
& = \left |\left( \frac{1}{i\omega \eta_0}\frac{\partial}{\partial r}  + 1 \right)  \delta u(\textbf{x}) \right |
\vspace{3mm}
\\
& \leq \frac{|\nabla A(\textbf{x}_0)| + O(r)}{\omega \eta_0} + \frac{|A(\textbf{x})|}{\eta_0} \left( r \left|  \widehat{\textbf{s}}^T \nabla^2 \phi ( \textbf{x}_0)  \widehat{\textbf{s}}  \right| + O(r^2)  \right)
\vspace{3mm}
\\
& + 2 |A(\textbf{x}_0)| \omega \left( \frac{1}{2} r^2 \left|  \widehat{\textbf{s}}^T \nabla^2 \phi ( \textbf{x}_0)  \widehat{\textbf{s}}  \right| + O(r^3)   \right)
+ 2r |\nabla A(\textbf{x}_0)| + O(r^2)
\vspace{3mm}
\\
& = \left( \frac{1}{\omega \eta_0} + 2r  \right) |\nabla A(\textbf{x}_0)|
+ \left( \frac{|A(\textbf{x})|r}{\eta_0} + |A(\textbf{x}_0)|\omega r^2  \right) \left|  \widehat{\textbf{s}}^T \nabla^2 \phi ( \textbf{x}_0)  \widehat{\textbf{s}}  \right|
\vspace{3mm}
\\
& + \frac{|A(\textbf{x})|}{\eta_0} O(r^2) + 2\omega |A(\textbf{x}_0)| O(r^3) + O(r^2).

\end{array}
\end{equation}

As shown in \ref{appendix-a}, on one hand $\delta U$ has to be small compared to $U$. On the other hand, the error in direction estimate from NMLA is $\cO(\frac{1}{wr})$. Assuming the smoothness of $A(\x)$ and $\phi(\x)$, i.e., boundedness of $\nabla A(\textbf{x})$, $A(\textbf{x})$ and $\nabla^2 \phi (\textbf{x})$, the leading term in $\delta U$ is $\omega r^2  |A(\textbf{x}_0)|\left|  \widehat{\textbf{s}}^T \nabla^2 \phi ( \textbf{x}_0)  \widehat{\textbf{s}}  \right|$ as $\omega \to \infty$, where $\widehat{\textbf{s}}^T \nabla^2 \phi ( \textbf{x}_0)  \widehat{\textbf{s}} $ is the curvature of the wave front.
Hence the radius of the sampling circle can at most be chosen  $r \sim \cO(\frac{1}{\sqrt{\omega}})$ as $\omega \to \infty$. Let
\begin{equation}
\begin{array}{ll}
& r = \frac{C_{\epsilon}}{\sqrt{\omega}},
\qquad |\nabla A(\textbf{x})| \leq C_1,
\qquad |A(\textbf{x})| \leq C_2,
\qquad \left| \textbf{s}^T \nabla^2 \phi(\textbf{x}) \textbf{s} \right| \leq C_3,

\end{array}
\end{equation}
Then
\begin{equation}
\begin{array}{ll}
|\delta U(\textbf{x})|
& \leq  2 C_{\epsilon}^2 C_3 \left |A(\textbf{x}_0) \right|
+ O \left( \frac{1}{\sqrt \omega} \right)
% =  2 C_{\epsilon}^2 C_3 \left |B_1 \right|
%+ O \left( \frac{1}{\sqrt \omega} \right)
\end{array}
\end{equation}

Choose $C_{\epsilon}$ small enough such that , say $2C_{\epsilon}^2 C_3 \leq \frac{1}{4}$, then the perturbation $\delta U(\textbf{x})$ satisfies the condition \ref{eq:NMLA_error_condition} for $\omega$ large enough, which implies  the error in ray direction estimate by NMLA is $\cO(\omega^{-\frac{1}{2}})$.
 %$\frac{\pi}{\eta_0 C_{\epsilon} \sqrt{\omega}}$ for large $\omega$.

\begin{remark}
The above analysis also shows that NMLA can not be used to estimate ray directions within a few wavelengths away from the point source since the curvature of the wave front is of order $\cO(w)$.
\end{remark}

\end{document}